\documentclass[a4paper, 12pt]{article}

\usepackage{ucs}                
\usepackage[utf8x]{inputenc} 	

\usepackage[T1]{fontenc}     
\usepackage{lmodern}         

\usepackage[text={160mm,250mm},centering]{geometry}

\usepackage{rcs}  \RCS $Revision: 1.14 $ \RCS $Date: 2019/05/31 22:38:57 $ \RCS $Author: hgm $ \RCS $RCSfile: 14_tensor_post-proc.tex,v $ \RCS $Id: 14_tensor_post-proc.tex,v 1.14 2019/05/31 22:38:57 hgm Exp $

\usepackage[sfdefault=cmbr,OMLmathsans]{isomath}
\usepackage{upgreek}
\usepackage{color}
\usepackage{amsmath}
\usepackage{amssymb}
\usepackage{amsbsy}
\usepackage{amsthm}
\usepackage{wasysym}
\usepackage{varioref}
\usepackage{tocbibind}          
\usepackage{float}
\usepackage{algpseudocode}
\usepackage{algorithm}
\usepackage{tensor}
\usepackage[breaklinks]{hyperref}
\usepackage{multirow} 
\usepackage{array}
\usepackage{graphicx}
\usepackage{caption}

\usepackage{authblk}
\usepackage{makeidx}
\makeindex

\usepackage[british]{babel}
\usepackage{refdef}
\usepackage{ifontdef}

\usepackage{amsmath}

\newcommand{\ignore}[1]{}

\newcommand{\Lp}{\mrm{L}}

\def\frks#1#2{{\raise0.8ex\hbox{{\leavevmode$\scriptstyle #1$}}}{\kern-0.4ex\hbox{/}}{\raise-0.2ex\hbox{\kern-0.4ex\hbox{{\leavevmode$\scriptstyle #2$}}}}}

\newcommand{\vrho}{\varrho}
\newcommand{\vepsilon}{\varepsilon}
\newcommand{\vsigma}{\varsigma}

\newcommand{\vek}[1]{\mathchoice{\displaystyle\boldsymbol{#1}}
{\textstyle\boldsymbol{#1}}{\scriptstyle\boldsymbol{#1}}
{\scriptscriptstyle\boldsymbol{#1}}}

\newcommand{\mat}[1]{\mathchoice{\displaystyle\mathbf{#1}}
{\textstyle\mathbf{#1}}{\scriptstyle\mathbf{#1}}
{\scriptscriptstyle\mathbf{#1}}}

\newcommand{\ops}[1]{\mathchoice{\displaystyle\mathsf{#1}}
{\textstyle\mathsf{#1}}{\scriptstyle\mathsf{#1}}
{\scriptscriptstyle\mathsf{#1}}}

\newcommand{\tnb}[1]{\mathchoice{\displaystyle\mathboldsans{#1}}
{\textstyle\mathboldsans{#1}}{\scriptstyle\mathboldsans{#1}}
{\scriptscriptstyle\mathboldsans{#1}}}

\newcommand{\tns}[1]{\mathchoice{\displaystyle\mathsans{#1}}
{\textstyle\mathsans{#1}}{\scriptstyle\mathsans{#1}}
{\scriptscriptstyle\mathsans{#1}}}

\newcommand{\vtil}[1]{\vek{\tilde{#1}}}
\newcommand{\vhat}[1]{\vek{\hat{#1}}}

\newcommand{\Tcek}[1]{\tnb{\check{#1}}}

\newcommand{\Ttil}[1]{\tnb{\tilde{#1}}}
\newcommand{\That}[1]{\tnb{\hat{#1}}}

\newcommand{\Shat}[1]{\tns{\hat{#1}}}

\newcommand{\EXP}[1]{\mathbb{E}\left(#1\right)}

\DeclareMathOperator{\diag}{diag}

\DeclareMathOperator{\spn}{span}

\DeclareMathOperator{\sign}{sign}

\newcommand{\dd}{\mathop{}\!\partial}
\newcommand{\di}{\mathop{}\!\mathrm{d}}

\newcommand{\citep}[1]{\cite{#1}}
\newcommand{\KL}{Karhunen-Lo\`eve}

\newcommand{\bkt}[2]{\langle #1 | #2 \rangle}

\newcommand{\ns}[1]{| #1 |}
\newcommand{\nd}[1]{\| #1 \|}

\newcommand{\trpos}{{\ops{T}}}
\definecolor{myred}{rgb}{1, 0.2, 0.2}

\newcommand{\R}{\mathbb{R}}
\newcommand{\N}{\mathbb{N}}

\newcommand    {\ten}    { \otimes }
\newcommand    {\Ten}    { \bigotimes }

\newcommand{\BIGOP}[1]{\mathop{\mathchoice%
{\raise-0.22em\hbox{\huge $#1$}} {\raise-0.05em\hbox{\Large $#1$}}
{\hbox{\large $#1$}}{#1}}}

\newcommand{\BIGboxplus}{\mathop{\mathchoice%
{\raise-0.35em\hbox{\huge $\boxplus$}}%
{\raise-0.15em\hbox{\Large $\boxplus$}}{\hbox{\large
$\boxplus$}}{\boxplus}}}
\newcommand{\bigtimes}{\BIGOP{\times}}

\newtheorem{theorem}{Theorem}[section]
\newtheorem{remark}[theorem]{Remark}

\newtheorem{definition}[theorem]{Definition}

\newtheorem{example}[theorem]{Example}

\newcommand{\thetitle}{Post-Processing of High-Dimensional Data}
\newcommand{\authorme}{Mike Espig}
\newcommand{\authorwh}{Wolfgang Hackbusch}
\newcommand{\authoral}{Alexander Litvinenko}
\newcommand{\authorhgm}{\authorcr Hermann~G. Matthies}
\newcommand{\authorez}{Elmar Zander}

\newcommand{\affilwire}{\small 
                        Technische Universit\"at Braunschweig, Brunswick, Germany}
\newcommand{\affilmis}{\small Max-Planck-Institute for Mathematics
                       in the Sciences (MIS), Leipzig, Germany}
\newcommand{\affilrwth}{\small Rheinisch-Westf\"alische Technische Hochschule (RWTH)
                         Aachen, Germany}
\newcommand{\affilzwi}{\small 
Wests\"achsische Hochschule Zwickau, 
Zwickau, Germany}
                       
\newcommand{\theauthor}{\authorme, \authorwh, \authoral, \authorhgm, \authorez}

\newcommand{\textdate}{\today}

\begin{document}

\title{\thetitle
     \thanks{Work partly supported by the Deutsche
         Forschungsgemeinschaft (DFG) and the Alexander von Humboldt Foundation (AvH).}
          }

\author[a]{\authorme}
\author[b]{\authorwh}
\makeatletter
\author[c]{\authoral
\thanks{Corresponding author: RWTH Aachen, 52072 Aachen,
       Germany,\\  e-mail: \texttt{Litvinenko@uq.rwth-aachen.de}}
}
\makeatother
\author[d]{\authorhgm}
\author[d]{\authorez}

\affil[a]{\affilzwi}
\affil[b]{\affilmis}
\affil[c]{\affilrwth}
\affil[d]{\affilwire}


\maketitle



\begin{abstract}
Scientific computations or measurements may result in huge volumes
of high-dimensional data, for instance $10^{20}$ or $100^{300}$ elements. Often these can be thought of representing a real-valued
function on a high-dimensional domain, and can be conceptually
arranged in the format of a tensor of high degree in some truncated or lossy
compressed format.  We look at some common post-processing tasks
which are too time and storage consuming in the uncompressed data format and not obvious in the compressed format, as such
huge data sets can not be stored in their
entirety, and the value of an element
is not readily accessible through simple look-up.
The tasks we consider are finding the location of maximum or minimum,
or minimum and maximum of a function of the data, or
finding the indices of all elements in some interval --- i.e.\ level
sets, the number of elements with a value in such a level set, the
probability of an element being in a particular level set, and the
mean and variance of the total collection.
The algorithms to be described are fixed point iterations of particular
functions of the tensor, which will then exhibit the desired result.
For this, the data is considered as an element of a high
degree tensor space, although in an abstract sense, the algorithms are
independent of the representation of the data as a tensor.  
All that we require is: the data are given in a compressed data format, and the rank truncation procedure preserves compression.
We allow the actual computational representation to be a lossy compression, and we
allow the algebra operations to be performed in an approximate fashion, so
as to maintain a high compression level.  One such example which we address explicitly
is the representation of data as a tensor with compression in the form of a low-rank representation.



\end{abstract}

%
%
%
%
%
%
%
%
%
%


\pagenumbering{roman}

\cleardoublepage

\tableofcontents
\cleardoublepage



\pagenumbering{arabic}

\section{Introduction}  \label{S:intro}
Many scientific and engineering computations or measurements, as well as economic
or financial applications, produce large volumes of data.  For simplicity, assume that
we have computed or observed large quantities of one real-valued variate.
Assume further that the amount of data is so large that it can not be held or
stored in its entirety and has to be compressed in some way. 

Frequently, the task arises for example to find the location with the
maximum value of the data set, or to find the locations of all values which lie
in a given interval --- we call this a level set.  Other, similar tasks we consider
are finding the number of values in a given interval, the probability of being in
a given interval, or finding the mean or the variance of the data set. All these tasks can be used, for example, to construct high-dimensional histograms and probability density functions.
These tasks above are trivial for small data sets, typically performed by inspecting each
datum.  But if we consider truly huge amounts of data which can not be stored
in full because this would exceed the capacity of any storage device, 
but only in a compressed manner, and where additionally it may not be
possible due to time constraints to inspect each element, these tasks
are not trivial any more.  This is due to the fact that
in the compressed representation the data values are normally not directly accessible
and require additional processing.  Such compression will
in general be ``lossy'', so that not each value can be restored exactly, but only
up to a certain reconstruction error.

Assuming the data as an element of some set $\C{T}$, 
the algorithms are independent of the representation of the data, as well as
from the compression and reconstruction
technique used, subject only to the possibility of approximately performing the
operations of an Euclidean associative commutative algebra.  This means that we assume
that $\C{T}$ is a vector space with an inner product, and additionally
an associative commutative bilinear multiplication, making it into an algebra.
The simplest example of such a structure is to envisage the data points as a vector
$\vek{w}\in\D{R}^N$ with an appropriate $N\in\D{N}$, reflecting the amount of data.
The vector space operations are clear, as is the canonical Euclidean product,
and for the commutative multiplication consider the point-wise or Hadamard product,
i.e.\ the component-wise multiplication of two vectors.  We shall allow that all these
operations are performed only approximately, in order to maintain a high compression level.

Large volumes of data, especially when they can be thought of as samples or discrete
values of some real-valued function on a high-dimensional space, can often be
arranged in form of a tensor \citep{Kolda:01, Kolda:07, DMV-SIAM2:00, HA12}.
This offers the possibility to use approximate or compressed tensor representations.
Here we will especially show one possibility, namely low-rank
representations, which generalise the truncated singular value decomposition for matrices.
For the sake of completeness, we shall show the possible implementation of the above
mentioned algebraic operations for some of the more common low-rank formats,
and the approximated, compressed, or truncated representation of their results;
see also \citep{Oseledets-S-T-2009, Espig2019}.

The proposed algorithms are iterative in nature,
and the convergence tolerance can be adapted to the reconstruction error.
The basic idea for the algorithms, which operate only on the algebraic structure, is
an iterative scheme which converges to a result which solves the desired problem
in some way.  Such iterations typically destroy the compressed representation,
so they have to be combined with re-compression or truncation.

\subsection{An example} \label{SS:xmpl-motiv}

In many chemical applications one has $n=100$, and $d$ is few hundreds or even 
thousands \citep{MikeRaoKhor}, making in total a tensor with $N=n^d$ entries. 
Some other high-dimensional problems from chemistry and physics like 
Hartree-Fock-, Schr\"odinger-, or Master-equations, and low-rank tensor methods 
of their solution are for example considered in 
\citep{KhorSurvey, khoromskaia2018tensor, khorBook18, Dolgov_Master}.
Another example of large volumes of high-dimensional data are  
satellite data.  Satellites collect data over a very large areas (e.g.\ the data collected 
by the National Center for Atmospheric Research (USA) \citep{GrossmanMazzucco2002}.
Big data can also come from a computer simulator codes such as a solution of a
multi-parametric equation, e.g.\ weather research and forecasting, and
climate models \citep{Eyring2019}.  Also
Oil$\&$Gas companies daily collect sensor
data from multiple sources.  
Another source for huge volumes of data are high-energy
particle accelerators like CERN \citep{Brumfiel2011}.

%
\ignore{For the sake of simplicity, we will treat all arguments of the function $v$ in
\refeq{eq:XXaX} equally, denoting them in that case
by}
Let $\vek{p}$ be a function, such that $\vek{p}=(x,t,\vek{\mu})=(p_1,\dots,p_d)$.
Wir consider a real-valued function $\vek{p}\mapsto w(\vek{p})=w(p_{1},\ldots,p_{d})$
and its samples
\begin{equation}  \label{eq:XXbX}
   \tnb{w} = (w(\vek{p}_{\vek{m}}))_{\vek{m}\in\D{N}^d} = 
   (w(p_{m_1},\ldots,p_{m_d}))_{\vek{m}} = (\tns{w}_{\vek{m}})_{\vek{m}} =
   (\tensor*{\tns{w}}{_{m_1}_{\dots}_{m_d}})_{m_1 \dots m_d} \in 
   \bigotimes_{\ell=1}^d \D{R}^{M_\ell},
\end{equation}
which form a tensor of order or degree $d$ representing a
total of $N:=\prod_{\ell=1}^d M_\ell$ values, and where for the sake of convenience
multi-indices $\vek{m}\in\C{M} := \bigtimes_{\ell=1}^d \{1,\dots,M_\ell\}\subset\D{N}^d$
have been introduced.

As one may see, a tensor can be simply defined as a high-order matrix or multi-index 
array, where multi-indices are used instead of indices.  As an
example, assume that we have $N=10^{16}$ data points.  Now a vector
$\vek{w}\in\D{R}^N$ is just a tensor of first order, but one may \emph{reshape}
the data, say into a matrix $\vek{W}\in\D{R}^{10^8 \times 10^8}\cong \D{R}^{10^8}\otimes\D{R}^{10^8}$
 --- a tensor of 2nd order.  Reshaping further, one may take
$\D{R}^{10^4 \times 10^4 \times 10^4 \times 10^4}\cong \bigotimes_{k=1}^4 \D{R}^{10^4}$
--- a tensor of 4th order --- as well as other combinations, such as 
$\tnb{w}\in\bigotimes_{k=1}^{16} \D{R}^{10}$ --- an array of size $10 \times \dots \times 10$ 
(16 times) --- which is a tensor of 16th order; and finally all the way to
$\tnb{w}\in\bigotimes_{k=1}^{16} (\D{R}^{2}\otimes\D{R}^{5})$ --- a tensor of order 32.

The tensors  obtained in this way contain not only rows and columns, but also
\emph{slices} and \emph{fibres} \citep{Kolda:01, Kolda:07, DMV-SIAM2:00, HA12}. 
These slices and fibres can be analysed for linear dependencies, super symmetry,
or sparsity,  and may result in a strong data compression.  To have a first glimpse
of possible compression techniques, assume that in the above example the data has
been stored in the matrix $\vek{W}\in\D{R}^{M\times M}$, where $M=10^8=\sqrt{N}$,
and consider its singular value decomposition (SVD)
$\vek{W}=\vek{U}\vek{\Sigma}\vek{V}^\trpos = \sum_{m=1}^M \vsigma_m \vek{u}_m \vek{v}_m^\trpos$,
where  $\vek{\Sigma}=\diag(\vsigma_1,\dots,\vsigma_M)$ is the diagonal matrix
of singular values $\vsigma_m\ge 0$, assumed arranged by decreasing value,
and $\vek{U}=[\vek{u}_1,\dots,\vek{u}_M]$,
$\vek{V}=[\vek{v}_1,\dots,\vek{v}_M]$ collect the right and left singular vectors.
Often there is a number $r\ll M$, such that for some small $\vepsilon>0$ one has
$\vsigma_m \le \vepsilon$ for all $m>r$. Then one can formulate a compressed or
truncated version 
\[
\vek{W}\approx\vek{W}_r =\vek{U}_r\vek{\Sigma}_r\vek{V}_r^\trpos
=\sum_{m=1}^r \vsigma_m \vek{u}_m \vek{v}_m^\trpos
=\sum_{m=1}^r  \vek{w}^{(1)}_m \otimes \vek{w}^{(2)}_m,
\] 
where $\vek{\Sigma}_r=\diag(\vsigma_1,\dots,\vsigma_r)$,
$\vek{U}_r=[\vek{u}_1,\dots,\vek{u}_r]$, $\vek{V}_r=[\vek{v}_1,\dots,\vek{v}_r]$,
and $\vek{w}^{(1)}_m = \sqrt{\vsigma_m} \vek{u}_m$, $\vek{w}^{(2)}_m = \sqrt{\vsigma_m} \vek{v}_m$.
This reduced SVD is a special case of \refeq{eq:CPI}.  For the sake of a concrete example,
assume that $r=100$.  One may observe that the
required storage has been reduced from $N=M^2=10^{16}$ to $2\times r\times M = 2\times 10^{10}$
for the $2\,r$ vectors $\{\vek{w}^{(1)}_m, \vek{w}^{(2)}_m\}_{m=1}^r$.

The set of possible objects we want to work with will be denoted by $\C{T}$.
We assume that this set is equipped with the operations of an associative and
commutative real algebra, and carries an inner product.  It is a well known result
\citep{segalKunze78} that such algebras --- modulo some technical details --- are
isomorphic via the Gel'fand representation to a function algebra --- 
to be more precise continuous real valued functions on a compact set.  
Under point-wise addition and multiplication by scalars, such functions
clearly form a vector space, and if one includes point-wise multiplication of functions,
they form an associative and commutative real algebra.  In our case this is simply
the function algebra $(\C{M}=\bigtimes_{\ell=1}^d \{1,\dots,M_\ell \} \to \D{R})$, 
which is obviously the same as $\bigotimes_{\ell=1}^d \D{R}^{M_\ell}$.
The advantage of this abstract algebraic formulation is not only that it
applies to any kind of data representation on which one may define these algebraic
operations, but that one may use algorithms \citep{NHigham} which have been developed
for the algebra of real $N\times N$ matrices $\F{gl}(\D{R},N)
= \D{R}^{N \times N} \cong \D{R}^N\otimes\D{R}^N$. 

In our concrete examples we work directly with the data represented as tensors
$\C{T}:= \bigotimes_{\ell=1}^d \D{R}^{M_\ell}$.  Even without the identification
with $(\C{M} \to \D{R})$, it is obvious that this is naturally a vector space.  
As also $\C{T}\cong \D{R}^N$, it is clear that $\C{T}$ can be equipped with 
the canonical Euclidean inner product from $\D{R}^N$.  The associative and 
commutative algebra product comes from the identification with the function 
algebra $(\C{M} \to \D{R})$, i.e.\ the point-wise product, which is also known
as the \emph{Hadamard} product \citep{HA12}.  Hence it is clear that for 
data $\vek{w}\in\D{R}^N$ it is not difficult at all to define the algebraic 
operations, but rather how to perform them when the data is in a compressed format.

To simplify later notation on the number of operations and amount of storage
needed, we will often make the assumption that $M_1 = \dots = M_d = n$, so that the
tensor in \refeq{eq:XXbX} represents $N=n^d$ values.
As already mentioned, our example of data compression is based on
low-rank approximations to elements in $\C{T}$.
Although low-rank tensor data formats and techniques are almost
unavoidable if working with large high-dimensional data sets,
we would like to stress that the algorithms presented here are
formulated purely in terms of the abstract algebraic structure,
and are thus independent of the particular representation.

Whereas a general element $\tnb{w}\in\C{T}$ hence has $N=n^d$ terms,
a compressed representation --- some low-rank versions of which
will be discussed in more detail in \refS{tensor-rep} --- will have
significantly fewer terms.  For example, the \emph{CP}-decomposition
(\emph{canonical polyadic}) representation, truncated to $r$ terms,
\begin{equation}  \label{eq:CPI}
   \tnb{w} \approx \tnb{w}_r = \sum_{i=1}^r   \bigotimes_{k=1}^d 
     \vek{w}_i^{(k)}; \;  \vek{w}_i^{(k)} \in \D{R}^{n},\; \text{has }\,
    r\times n \times d \; \text{terms in}\; \tnb{w}_r.
\end{equation}
If the \emph{rank} $r$ is reasonably small compared to $N=n^d$ independent of $n$ and $d$,
then we have an approximation $\tnb{w}_r\in\C{T}$ with much less storage, which also
only depends on the \emph{product} of $n$ and $d$ and is not exponential in the
dimension $d$.  But now the maximum can not
be easily looked up; to get a particular element, the expression \refeq{eq:CPI}
has to be evaluated for that index combination.  In the following
we shall assume that we work with approximations such as in \refeq{eq:CPI}
which need much less storage.  One has to make sure then that the \emph{algebraic}
operations of the Hadamard algebra structure on $\C{T}$ can be performed efficiently
in an approximative manner, so that they
have much lower complexity than the obvious $\C{O}(n^d)$, although
the single elements of $\tnb{w}$ which are usually needed for point-wise operations
--- and post-processing, see \refSS{post-proc-intro} --- are not directly available.
Thus the motivating factors for applying compression, and in particular
low-rank tensor techniques, include the following:
\begin{itemize}
\item The storage cost is reduced, depending on the tensor format, from 
      $\C{O}(n^d)$ to $\C{O}(d r n)$, where $d>1$. 
\item The low-rank tensor approximation is relatively new, but already a well-studied technique 
      with free software libraries available. Other data  compression techniques are available.
\item The approximation accuracy is fully controlled by the tensor rank.
      The full rank gives an exact representation.
\item Even more complicated operations like the Fourier transform can be performed
      efficiently in low-rank format.  The basic fast Fourier transform on $\tnb{w}\in\C{T}$
      would have complexity $\C{O}(n^d \log (n^d))$.  These low-rank techniques can
      either be combined with the fast Fourier transform giving a complexity of $\C{O}(d r n\log n)$,
      or can even be further accelerated at the price of an additional approximation error
      yielding a \emph{superfast} Fourier transform  \citep{nowak2013kriging, DoKhSa-qtt_fft-2011}.
\end{itemize}
On the other hand, general limitations of a compression technique are that 
\begin{itemize}
\item it could be time consuming to compute a compression, or in particular a
      low-rank tensor decomposition; 
\item it requires an axes-parallel mesh; 
\item although many functions have a low-rank representation, in practice
      only some theoretical estimates exist.
\end{itemize}
But the fact still remains that there are situations where storage of
all items is not feasible, and some kind of compression has to be employed.

\subsection{Post-processing tasks} \label{SS:post-proc-intro}
This work is about exploiting the compression together with the structure of
a commutative algebra and the ability to perform the algebraic operations
in the compressed format, at least approximately.  The tensor product structure allows efficient calculations
to be performed on a sample of function like \refeq{eq:XXbX}
with the help of expressions such as \refeq{eq:CPI} or other compression
techniques.

This tensor product structure---in this case
multiple tensor product structure---is typical for such parametric
problems \citep{hgmRO-1-2018}.  What is often desired, is a representation which allows for
the approximate evaluation of the state of \refeq{eq:XXbX} as function of $\vek{\mu}\in \C{M}$
without actually solving the system again.  Furthermore, one would like this
representation to be inexpensive to evaluate, and for it to be
convenient for certain post-processing tasks, for example like finding
the minimum or maximum value over some or all parameter values.

The tasks we consider here are
\begin{enumerate}  
\item the indices or values of \emph{maxima, minima} of $\tnb{w}$,
\item the indices or values of \emph{maxima, minima} of some function $f$ of $\tnb{w}$,
\item the index or value of the datum of $\tnb{w}$ \emph{closest} to a given number,
\item level sets, i.e.\ (the indices of) all $\tensor*{\tns{w}}{_{m_1}_{\dots}_{m_d}}
      \in[\omega_0, \omega_1]$,
\item the \emph{number of indices} such that $\tensor*{\tns{w}}{_{m_1}_{\dots}_{m_d}}
      \in[\omega_0, \omega_1]$,
\item the mean (sum of all items) and variance (sum of all items squared),
\item some auxiliary functions, such as the algebraic inverse $\tnb{w}^{\odot-1}$ of 
      $\tnb{w}$ --- in our case the Hadamard inverse --- and others 
      like the $\sign(\tnb{w})$ function which will be seen to
      be needed for the above tasks in \refSS{post-proc-tasks}.
\end{enumerate}

%
As an example, consider finding the maximum value.
A naïve approach to compute the maximum would be to visit and inspect each element,
but then the number of visits is $\C{O}(n^d)$, exponential in $d$,
which may be unacceptable for high $d$.
Such phenomena when the total complexity/storage cost depends exponentially on the
problem dimension have been called the curse of dimensionality \citep{HA12}.

The idea for such iterative post-processing algorithms, 
in the form of finding the maximum, was first presented in \citep{EspigDiss}
for low-rank tensor approximations.  Some further post-processing tasks for such low-rank
tensor representations were investigated in \citep{ESHALIMA11_1}.  To give an example
from \citep{EspigDiss} of the possible savings in space and time, 
assume that one looks at the solution of a $d$-dimensional Poisson equation:
\[
  - \nabla^2 u = f \; \text{ on } \; \C{G}=[0,1]^d \; \text{ with } \; 
  u\vert_{\dd \C{G}} = 1,
\]
and right-hand-side 
\[
  f(x_1,\dots,x_d) \propto   \sum_{k=1}^d \prod_{\ell=1, \ell\ne k}^d  x_\ell (1-x_\ell) .
\]
Assume further that this is solved numerically by a  standard finite-difference method with
$n=100$ grid-points in each direction, so that the solution vector $\vek{u}$ has $N=n^d$ data
points in $\D{R}^N$, but can also naturally be viewed as a tensor $\tnb{u}$ of degree $d$ in 
$\bigotimes_{k=1}^d \D{R}^{100}$.  So for full storage, one needs $N=n^d$
storage locations, and if one is looking for the maximum by inspecting each element,
one would have to inspect $N$ elements.  

We mention in passing that of course
actually solving the discrete equation is a challenge for $d>3$,
but that is a different story, and for this we refer to \citep{EspigDiss}
and the references therein.

\begin{table}[h!]
\begin{center}
\caption{Computing times (4th column) on 2~GHz dual-core CPU to find maximum.}
\label{tab:find_max}
\begin{tabular}{r|l|l|l}
$d$ & \# loc's.:  & $\approx$ years [a]   & actual \citep{EspigDiss} time [s]\\
    &   $N=n^d$   & inspect. $N$ & of Algorithm~\ref{alg:power-it} \\
\hline
 $25$ & $10^{50}$  & $ 1.6\times 10^{33}$ & 0.16 \\
 $50$ & $10^{100}$ & $ 1.6\times 10^{83}$ & 0.42 \\
 $75$ & $10^{150}$ & $ 1.6\times 10^{133}$ & 1.16 \\
$100$ & $10^{200}$ & $ 1.6\times 10^{183}$ & 2.58 \\
$125$ & $10^{250}$ & $ 1.6\times 10^{233}$ & 4.97 \\
$150$ & $10^{300}$ & $ 1.6\times 10^{283}$ & 8.56
\end{tabular}

\end{center}
\end{table}%

Now assume for the sake of simplicity that
one can inspect $2\times 10^9$ elements per second ---
on an ideal 2~GHz CPU with one inspection per cycle.
Then for $\tnb{u} \in (\D{R}^n)^{\otimes d} \cong \D{R}^{n^d}$ 
the times needed to find the maximum for the full data-set
are shown in \reftab{tab:find_max} in the third column---assuming that it were somehow 
possible to store all the values indicated in the second column---whereas the
actual computation with a compressed format is shown in the last and fourth column.

It is obvious that for growing $d$ for the full representation
the computational complexity and the
storage requirements quickly become not only unacceptable, but totally
impossible to satisfy.  The second and third column behave like $\C{O}(n^d)$ and
grow exponentially with $d$, whereas for a low-rank representation 
$\tnb{u}_r$ of $\tnb{u}$ not only
can the data be stored on a modest laptop, but for the simple Algorithm~\ref{alg:power-it}
--- to be explained later in \refSS{auxil} --- the computing times in the fourth column 
are in terms of seconds and behave like $\C{O}(n\, d^3)$.

%
%

\subsection{State of the art} \label{SS:state-art}
This is not a discussion of the state of the art regarding
general tensor formats and their low-rank approximations, in quantum physics
also known as \emph{tensor networks}.  For the physical motivations and numerical
developments from there see \citep{VI03, Sachdev2010-a, EvenblyVidal2011-a, Orus2014-a, 
BridgemanChubb2017-a, BiamonteBergholm2017}.  
For the mathematical and numerical view we refer to 
the review \citep{Kolda:07}, the monographs \citep{HA12, khoromskaia2018tensor, khorBook18},
and to the literature survey on low-rank approximations \citep{GraKresTo:13}.
In the following, we concentrate on the question of post-processing such data.

The idea of finding the largest element and the corresponding location of
a tensor by solving an eigenvalue problem, where the matrix and the vectors are
tensors given in the CP tensor format, was introduced in \citep{EspigDiss}.  
Additionally, the first numerical schemes to compute the 
point-wise inverse and the sign function were introduced, as well as
the rank-truncation procedure for tensors given in the CP format.
Later, in \citep{ESHALIMA11_1}, these ideas were extended and applied to tensors
which were obtained 
after discretisation and solution of elliptic PDEs with uncertain coefficients. Having the sign function computed, we can compute (high-dimensional) histograms, which are needed for approximating high-dimensional probability density functions, cumulative distribution functions, quntiles and exceedance probabilities.
Another group of authors, in \citep{DoKhSavOs_mEIG:13}, by
combining the advances of the density matrix renormalisation 
group and the variational numerical renormalisation group methods,
approximated several low-lying
eigenpairs of large Hermitian matrices simultaneously in the block version of the 
tensor train (TT) format via the alternating minimisation of the block Rayleigh quotient 
sequentially for all TT cores.

In \citep{DolgLitv15, dolgov2014computation}, the authors suggested methods to
compute the mean, the variance, and sensitivity indices in the TT format 
with applications to stochastic PDEs, whereas the diagonalisation of large 
Toeplitz or circulant matrices via combination of the 
fast Fourier and the CP tensor format was shown in \citep{nowak2013kriging}.

An investigation of approximations to eigenfunctions of a certain class of elliptic 
operators in $\D{R}^d$ by finite sums of products of functions with separated 
variables is the topic of \citep{HKST_Eigs_12}.  Various tensor formats were used for 
a new class of rank-truncated iterative eigensolvers.  The authors were able to reduce 
the computational cost from $\C{O}(n^d)$ to $\C{O}(n)$.  They solved large-scale spectral problems  from quantum 
chemistry: the Schrödinger, the Hartree–Fock, and the Kohn–Sham equations in electronic 
structure calculations.

\subsection{Outline of the paper} \label{SS:outline}
In the following \refS{algos} the necessary material for abstract algebras
is quickly reviewed.  Then the algorithms and functions on the algebra $\C{T}$
used to compute the post-processing tasks outlined in \refSS{post-proc-intro}
are formulated in an abstract fashion, independent of the representation
chosen for the data.  Even as the formulation uses only the abstract algebra operations,
we point out and motivate what this means in terms of the Hadamard algebra.
But also for the Hadamard algebra, the presentation is independent of the tensor
format to be chosen.  Furthermore, for the iterative algorithm --- a fixed point
iteration --- it is discussed how the possible truncation operation influences
the convergence.

In \refS{tensor-rep} we present some concrete examples of compression of high-dimensional
data $\vek{w}\in\D{R}^N$, once it is identified with a tensor $\tnb{w}\in
\C{T}:= \bigotimes_{\ell=1}^d \D{R}^{M_\ell}$ with $N=\prod_{\ell=1}^d M_\ell$.
For such tensors of degree $d$ we discuss the \emph{canonical polyadic} (CP) representation, the \emph{tensor-train} (TT) representation (in appendix),
and the Tucker tensor formats \cite{Espig2019}.
In particular, we show how the algebra operations can be carried out in the
different tensor formats, the numerical effort involved, and the effect these
algebraic operations have on the compression or low-rank representation.
As the compression level may deteriorate, i.e.\ the rank of the approximation
may grow, it is important that one is able to re-compress.  Pointers to such
methods in the literature are included as well.

The \refS{num-exp} contains a few numerical examples and a discussion of their
results used to illustrate the algorithms
of \refS{algos} and the representations from \refS{tensor-rep}.  These examples
are solutions of elliptic high-dimensional or parametric resp.\ stochastic partial
differential equations, a domain where the data is naturally in a tensor format.
But as already pointed out, for any data it is just a --- often only conceptual ---
\emph{reshape} to consider it as an element of a tensor space.  The conclusion 
is contained in \refS{concl}.

\section{Algorithms for post-processing} \label{S:algos}
Given a data-set $\tnb{w}$ in some compressed tensor format,
where looking up each element is not trivial, and where additionally
it may not be computationally feasible to look at each element,
we still want to be able to perform tasks which essentially require looking at
each datum.  The \emph{tasks} 
we consider are outlined in \refSS{post-proc-intro}.
These tasks and auxiliary functions will be computed by finding an
\emph{iteration} (mapping) $\tns{\Phi}_P$ for each post-processing task $P$,
such that the \emph{fixed-point} $\tnb{v}_*$ of $\tns{\Phi}_P$ 
--- i.e.\ $\tns{\Phi}_P(\tnb{v}_*)=\tnb{v}_*$ --- is either the
sought solution for the task, or computes one of the auxiliary functions.

\subsection{Iteration with Truncation} \label{SS:itrunc}
When performing computations using the operations of the algebra --- 
like the Hadamard algebraic operations on tensors in some compressed format --- 
after the operation the compression will typically
be sub-optimal, which means that the result has to be compressed again.
Thus, when performing the algebraic operations for a fixed-point iteration, the compression
would either get worse and worse in each iteration, or one has to re-compress
or \emph{truncate} the result again to a good compression level.
Therefore a re-compression after each or after a number of algebraic operations
may be necessary, and we thus allow that the algebraic operations are possibly
only executed approximately.

In our example case the low-rank representation of tensors explained
in \refS{tensor-rep} acts as compression, and the rank may increase
through the algebraic operations, and hence we will use \emph{truncation}
$\tns{T}_\epsilon$ to \emph{low rank} $r$ with error $\epsilon$
\citep{HA12, Espig:2013} to avoid the problem of a deteriorating
compression level.

In other words, with a general compressed representation $\tnb{w}_r$ of $\tnb{w}$
the computation will be a \emph{truncated} or \emph{perturbed iteration} \citep{whBKhEET:2008,
Zander10}.  If we denote the general compression mapping by $\tns{T}_\epsilon$ ---
meaning compression with an accuracy $\epsilon$ --- the iteration map
is changed from $\tns{\Phi}_P$ to $\tns{T}_\epsilon \circ \tns{\Phi}_P$.

The general structure of the iterative
algorithms for a post-processing task $P$ 
or for an auxiliary function is shown in Algorithm~\ref{alg:basic}.
  Here we want to collect
some results for this kind of iteration.
\begin{algorithm}
  \caption{Iteration with truncation}
       \label{alg:basic}
   \begin{algorithmic}[1]    
    \State Start with some initial \emph{compressed} guess  $\tnb{v}_0$
           depending on task $P$.
    \State $i\gets 0$
    \While  {\emph{no convergence}} 
       \State $\tnb{z}_{i} \gets \tns{\Phi}_P(\tnb{v}_i)$;\Comment{See Remark below}
       \Statex\Comment{ the iterator $\tns{\Phi}_P$ may \emph{deteriorate} the
                             compression level } 
       \If {compression level of $\tnb{z}_i$ is too bad} 
          \State $\tnb{v}_{i+1} \gets \tns{T}_\epsilon(\tnb{z}_i)$;
          \Statex\Comment{ use \emph{truncation} $\tns{T}_\epsilon$ to
                                  compress $\tnb{z}_i$ with error $\epsilon$ }
       \Else
          \State $\tnb{v}_{i+1} \gets \tnb{z}_i$;  
       \EndIf
       \State $i\gets i+1$
    \EndWhile 
  \end{algorithmic}
\end{algorithm}
For such a truncated or perturbed iteration as in  Algorithm~\ref{alg:basic},
it is known that
\begin{enumerate}
\item if the iteration by $\tns{\Phi}_P$ is \emph{super-linearly}
   convergent, the \emph{truncated iteration}  $\tns{T}_\epsilon \circ \tns{\Phi}_P$ will
   still \emph{converge} super-linearly, but finally \emph{stagnate} in an 
   $\epsilon$-neighbourhood of the fixed point $\tnb{v}_*$ \citep{whBKhEET:2008}.
   One could loosely say that the super-linear convergence of $\tns{\Phi}_P$ is
   stronger than the truncation by $\tns{T}_\epsilon$. 
\item if the iteration by $\tns{\Phi}_P$ is \emph{linearly} convergent with \emph{contraction}
   factor $q$, the \emph{truncated iteration}  $\tns{T}_\epsilon \circ \tns{\Phi}_P$ will
   still \emph{converge} linearly, but finally
   \emph{stagnate} in an $\epsilon/(1-q)$-neighbourhood
   of $\tnb{v}_*$ \citep{Zander10}.  Again, one could loosely say that
   iteration by $\tns{\Phi}_P$ and truncation by $\tns{T}_\epsilon$ balance each other,
   thus resulting in a larger neighbourhood of stagnation. 
\end{enumerate}

We shall assume that the truncation level has thus been chosen according to the
desired re-construction accuracy and taking into account the possible influence
due to the convergence behaviour of the iterator $\tns{\Phi}_P$.
\begin{remark}
In Algorithm~\ref{alg:basic} we assume that the tensor rank of $\tnb{z}_{i}$ is not increasing too much after applying 
the mapping $\tns{\Phi}_P$ to the tensor $\tnb{v}_i$. This could be a strong assumption. Additionally, we assume that there is a 
mapping $\tns{T}_\epsilon$, which is able to truncate the rank of $\tnb{z}_i$.
\end{remark}
\subsection{Preliminaries and basic algebraic operations} \label{SS:prelim}
As already pointed out, we assume in general that the set $\C{T}$ is a vector
space.   Our example of the space $\C{T}=\bigotimes_{\ell=1}^d \D{R}^{M_\ell}$
of tensors of interest was already introduced in
\refSS{xmpl-motiv}.  This is clearly a vector space, so we are
able to add two such tensors, and multiply each by some real number, in other words
we may form linear combinations.  The additive neutral element --- the zero
tensor --- will be denoted as $\tnb{0}$.  
It is important that the compressed storage format
chosen allows to perform the vector space operations without going into the \emph{full}
representation.  This will be shown for some of the familiar tensor formats in
\refS{tensor-rep}.  Furthermore, as
$\C{T} := \bigotimes_{\ell=1}^d \D{R}^{M_\ell} \cong \D{R}^{M_1 \times \dots \times M_d}
 \cong \D{R}^N$ as vector spaces, we can carry the canonical Euclidean inner product on 
$\D{R}^N$ to $\C{T}$, which for $\tnb{u}, \tnb{v}\in\C{T}$ is denoted by 
\[
\bkt{\tnb{u}}{\tnb{v}}_{\C{T}} := \sum_{m_1=1,\dots,m_d=1}^{M_1,\dots,M_d} 
     \tns{u}_{m_1,\dots,m_d} \cdot \tns{v}_{m_1,\dots,m_d}
     = \sum_{\vek{m}\in\C{M}} \tns{u}_{\vek{m}}\cdot\tns{v}_{\vek{m}} .
\]
This makes $\C{T}$ into a Euclidean or Hilbert space.  We assume that the
computation of this inner product is
feasible when both $\tnb{u}$ and $\tnb{v}$ are in a compressed representation; again
this will be demonstrated in \refS{tensor-rep} for the low-rank representations
we shall consider as examples.

On this space $\C{T}$ an additional
structure is needed, namely a \emph{multiplication} which will make
it into a unital associative and commutative algebra.  In our example this will be
the \emph{Hadamard} multiplication, which apparently goes
back to \emph{Schur}, which is the point- or component-wise multiplication on 
$\bigotimes_{\ell=1}^d \D{R}^{M_\ell}\cong \D{R}^{\bigtimes_{\ell=1}^d M_\ell}$.
It is the usual point-wise multiplication of functions, in this case
the functions $\bigtimes_{\ell=1}^d M_\ell \to \D{R}$.
For $\tnb{u}, \tnb{v} \in \bigotimes_{\ell=1}^d \D{R}^{M_\ell}$, 
define the Hadamard product as
\begin{equation}  \label{eq:Hadamard-def}
  \odot: \C{T}\times\C{T}\to\C{T}, \quad \tnb{u}\odot\tnb{v} \mapsto \tnb{w} =
  (\tensor*{\tns{w}}{_{m_1}_{\dots}_{m_d}}) := 
  (\tensor*{\tns{u}}{_{m_1}_{\dots}_{m_d}}\cdot\tensor*{\tns{v}}{_{m_1}_{\dots}_{m_d}}).
\end{equation}
It is a bi-linear operation (linear in each entry), and it is commutative.
Thus this product makes $\C{T}$ into an \emph{associative} and \emph{commutative}
\emph{algebra}.  The symbol $\odot$ will be used both for the abstract
algebra product on $\C{T}$, as well as for the Hadamard product
on $\bigotimes_{\ell=1}^d \D{R}^{M_\ell}$.
As is any such algebra, for any $\tnb{w}\in\C{T}$ it holds that $\tnb{w}\odot\tnb{0}=\tnb{0}$.

Especially, as a function of $\tnb{v}$, the product $\tnb{u}\odot\tnb{v}$
is a linear map $\tnb{L}_{\tnb{u}} \in \E{L}(\C{T})$ on $\C{T}$.  
This is the familiar canonical representation
$\C{T}\ni\tnb{u} \mapsto \tnb{L}_{\tnb{u}} \in \E{L}(\C{T})$
of an associative algebra as linear maps on itself, i.e.\ in this case in a commutative
sub-algebra of the algebra of all linear maps $\E{L}(\C{T})$ with concatenation 
as product.

It easy to see that this Hadamard algebra has a unit element 
$\tnb{1}\in\bigotimes_{\ell=1}^d \D{R}^{M_\ell}$ for the product, 
\begin{equation}  \label{eq:Hadamard-unit}
  \tnb{1} := (1_{{m_1{\dots}m_d}}), \quad \forall \tnb{w}\in\C{T}:\; 
  \tnb{w} = \tnb{1}\odot\tnb{w} = \tnb{w}\odot\tnb{1},
\end{equation}
 a tensor with all entries equal
to unity --- this makes $\bigotimes_{\ell=1}^d \D{R}^{M_\ell}$ into a \emph{unital} algebra.
Observe that by defining for all $1\leq k \leq d$ the \emph{all-ones} vectors
$\vek{1}_{M_k} = [1,\dots,1]^T\in\D{R}^{M_k}$, the Hadamard unit has rank \emph{one}
according to \refeq{eq:CPI}: $\tnb{1} = \bigotimes_{k=1}^d \vek{1}_{M_k}$.
Obviously one has $\tnb{L}_{\tnb{1}} = \tnb{I}_{\C{T}}$, the identity in $\E{L}(\C{T})$.

Having defined a unit or neutral element for multiplication, we say that
$\tnb{w}\in\C{T}$ has a \emph{multiplicative inverse} iff there is an element,
denoted by $\tnb{w}^{\odot-1}\in\C{T}$, such that
\begin{equation}  \label{eq:Hadamard-inverse}
  \tnb{1} =  \tnb{w}^{\odot-1}\odot\tnb{w} = \tnb{w}\odot\tnb{w}^{\odot-1}.
\end{equation}
Obviously not all elements have an inverse, e.g.\ the zero element $\tnb{0}\in\C{T}$
is never invertible.
In any case, not every $\tnb{w}\in\bigotimes_{\ell=1}^d \D{R}^{M_\ell}$ has a
Hadamard inverse, for this
it is necessary that all entries $\tensor*{\tns{w}}{_{m_1}_{\dots}_{m_d}}\neq 0$
are non-zero, and then $\tnb{w}^{\odot-1} = 
(\tensor*{\tns{w}}{_{m_1}_{\dots}_{m_d}^{-1}}) = 
(\tensor*{(1/\tns{w})}{_{m_1}_{\dots}_{m_d}})$.  Note that 
$(\tnb{w}^{\odot-1})^{\odot-1}=\tnb{w}$ and $\tnb{1}^{\odot-1}=\tnb{1}$
as it should be, and it is easily seen that 
$\tnb{L}_{\tnb{w}^{\odot-1}} = \tnb{L}_{\tnb{w}}^{-1}$,
the inverse of $\tnb{L}_{\tnb{w}}$ in $\E{L}(\C{T})$.

Elements of the algebra which can be written as a square 
$\tnb{w}=\tnb{u}\odot\tnb{u} = \tnb{u}^{\odot 2}$ are called positive, they form
a convex cone $\C{T}_+ \subset \C{T}$; note that obviously the zero tensor $\tnb{0}$
and the multiplicative unit $\tnb{1}$ are positive, and that if an invertible
$\tnb{w}$ is positive, so is its inverse $\tnb{w}^{\odot -1}$.  As usual,
this defines an order relation $\tnb{u} \le \tnb{v} \Leftrightarrow \tnb{v}-\tnb{u}\in\C{T}_+$,
which implies that for any positive $\tnb{w}\in\C{T}_+$ one has $\tnb{0}\le\tnb{w}$.
Further one may observe that $\tnb{u}, \tnb{v}\in\C{T}_+ \Rightarrow
\tnb{u}\odot\tnb{v}\in\C{T}_+$.

A certain compatibility of the inner product and the algebra product is needed,
as we require that 
\begin{equation}   \label{eq:self-adj-mult}
\bkt{\tnb{w}\odot\tnb{u}}{\tnb{v}}_{\C{T}} = \bkt{\tnb{L}_{\tnb{w}}\tnb{u}}{\tnb{v}}_{\C{T}} =
\bkt{\tnb{u}}{\tnb{L}_{\tnb{w}}\tnb{v}}_{\C{T}} = \bkt{\tnb{u}}{\tnb{w}\odot\tnb{v}}_{\C{T}},
\end{equation}
i.e.\ that the action of the algebra product is \emph{self-adjoint}, or, in other words,
that the maps $\tnb{L}_{\tnb{w}}$ are self adjoint.  This condition is satisfied
for the Hadamard algebra.

The inner product with the unit element $\tns{\phi}(\tnb{w}) := \bkt{\tnb{w}}{\tnb{1}}_{\C{T}}$
defines, as function of $\tnb{w}$, a positive linear functional $\tns{\phi}$,
as for $\tnb{w}\in\C{T}_+$ one has
\[
  \tns{\phi}(\tnb{w}) = \tns{\phi}(\tnb{u}\odot\tnb{u})=\bkt{\tnb{u}\odot\tnb{u}}{\tnb{1}}_{\C{T}}
  = \bkt{\tnb{u}}{\tnb{u}}_{\C{T}} \ge 0 .
\]
It is a kind of ``trace'' or un-normalised state functional on the algebra, 
as $\tns{\phi}(\tnb{1}) = N = \dim \C{T}$ in the Hadamard algebra.  In particular
\begin{equation}   \label{eq:state-inn-pr}
\tns{\phi}(\tnb{u}\odot\tnb{v}) = \bkt{\tnb{u}\odot\tnb{v}}{\tnb{1}}_{\C{T}}=
  \bkt{\tnb{u}}{\tnb{v}}_{\C{T}}.
\end{equation}
Conversely, if the algebra comes equipped with such an un-normalised positive state functional,
the \refeq{eq:state-inn-pr} may be taken as the definition of an inner product,
which automatically satisfies \refeq{eq:self-adj-mult}.

We shall assume that we can compute the algebraic operations 
--- at least approximately --- in compressed representation (this will be shown
for the Hadamard algebra 
for the low-rank formats in \refS{tensor-rep}), and that we can compute the
multiplicative inverse also in compressed representation --- again at least approximately.

As in any unital algebra, if for $\lambda\in\D{C}$ the element $\tnb{w}-\lambda\tnb{1}$ fails
to be invertible --- this means that $\tnb{L}_{\tnb{w}} - \lambda\tnb{I}_{\C{T}}$
is not invertible in $\E{L}(\C{T})$ --- we shall say that $\lambda$ is in the
\emph{spectrum} of $\tnb{w}$; here --- as we are in finite dimensional
spaces --- it is an \emph{eigenvalue}, i.e.\ there is an
\emph{eigenvector} $\tnb{v}_\lambda\in\C{T}$
such that $\tnb{L}_{\tnb{w}} \tnb{v}_\lambda = \tnb{w} \odot \tnb{v}_\lambda = 
\lambda \, \tnb{v}_\lambda = \lambda\, \tnb{I}_{\C{T}} \tnb{v}_\lambda$.
From this discussion it is immediate that in the Hadamard algebra case, if $\lambda$ is an
eigenvalue of $\tnb{w}=(\tns{w_{\vek{m}}})\in\C{T}$ or $\tnb{L}_{\tnb{w}}\in\E{L}(\C{T})$,
there must be a multi-index $\vek{m}_\lambda\in\C{M}$ such
that $\tns{w}_{\vek{m}_\lambda} = \lambda$, and hence in the real
Hadamard algebra $\C{T}$ one has that always $\lambda\in\D{R}$.
Thus the spectrum of $\tnb{w}$ resp.\ $\tnb{L}_{\tnb{w}}$ is $\sigma(\tnb{w}) =
\sigma(\tnb{L}_{\tnb{w}}) = \{ \tns{w_{\vek{m}}} \mid \vek{m}\in\C{M} \}$.
Observe that according to \refeq{eq:self-adj-mult}, each
representing map $\tnb{L}_{\tnb{w}}$ from above is self-adjoint,
which means that also in general all the spectra are real --- $\sigma(\tnb{w}) \subset \D{R}$.

Specifically, the Euclidean Hadamard algebra
$\C{T}=\bigotimes_{\ell=1}^d \D{R}^{M_\ell}$ which we have constructed
is obviously isomorphic --- as a Euclidean algebra --- to $\D{R}^N$ with the
canonical Euclidean inner product when
equipped with the point- or element-wise Hadamard product $\odot_N$.  Let us denote
this isomorphism by $\tns{V}:\C{T}\to\D{R}^N$; it implies some ordering
of the terms of each $\tnb{w}\in\C{T}$.

More enlightening and less obvious may
be the unital algebra isomorphism with the commutative sub-algebra of
\emph{diagonal} matrices $\F{diag}(\D{R},N)$ in the full matrix algebra $\F{gl}(\D{R},N)
= \D{R}^{N \times N} \cong \D{R}^N\otimes\D{R}^N \cong \E{L}(\D{R}^N)$
with the usual matrix multiplication.
Let us denote this isomorphism by $\tns{M}:\C{T}\to\F{diag}(\D{R},N)$.
If $\vek{w}=\tns{V}(\tnb{w}) \in\D{R}^N$ is the vector containing all the elements of the
tensor $\tnb{w}$, then $\tns{M}(\tnb{w}):=\vek{W}=\diag(\vek{w})=\diag(\tns{V}(\tnb{w}))\in
\F{diag}(\D{R},N) \subset \F{gl}(\D{R},N)$
is the corresponding diagonal matrix, i.e.\ $\tns{M}=\diag\circ \tns{V}$.  Hence, by choosing the
canonical Euclidean basis in $\D{R}^N$ and its image by $\tns{V}^{-1}$ as a basis in $\C{T}$,
the canonical representation $\tnb{L}_{\tnb{w}}\in\E{L}(\C{T})$ is itself represented
by the matrix $\vek{W}=\diag(\vek{w})=\tns{M}(\tnb{w}) \in \F{diag}(\D{R},N)
\subset \F{gl}(\D{R},N) \cong \E{L}(\D{R}^N)$.

Let $\tnb{v}\in\C{T}$ be another tensor, with associated $\tns{V}(\tnb{v})=
\vek{v}\in\D{R}^N$ and $\tns{M}(\tnb{v})=\vek{V}=\diag(\tns{V}(\tnb{v}))=
\diag(\vek{v})\in\F{diag}(\D{R},N)$, then
from the definition of the isomorphism $\tns{M}$:
\begin{multline}   \label{eq:H-alg-matrix}
  \tns{M}(\tnb{w}\odot\tnb{v})=\tns{M}(\tnb{w})\tns{M}(\tnb{v}) = \vek{W}\vek{V} = 
  \diag(\tns{V}(\tnb{w}))\diag(\tns{V}(\tnb{v}))\\
  =\diag(\vek{w})\diag(\vek{v}) = \diag\left(\tns{M}(\tnb{w}) \tns{V}(\tnb{v})\right)
 = \diag\left( \vek{W} \vek{v}\right) = \diag\left( \diag(\vek{w}) \vek{v}\right),
\end{multline}
which contains the important equality 
\begin{multline}  \label{eq:H-alg-matrix-2}
\tnb{w}\odot\tnb{v} = \tns{M}^{-1}(\vek{W}\vek{V}) = 
   \tns{M}^{-1}(\diag\left( \vek{W} \vek{v}\right))= \tns{M}^{-1}(\diag(\vek{w})\diag(\vek{v}))=\\
   \tns{M}^{-1}(\diag\left( \tns{M}(\tnb{w}) \tns{V}(\tnb{v})\right)) = 
   \tns{V}^{-1}(\vek{W} \vek{v})= \tnb{L}_{\tnb{w}} \tnb{v},
\end{multline}
from where one sees that 
\begin{equation}   \label{eq:H-alg-matrix-3}
\tns{V}(\tnb{w}\odot\tnb{v}) =  \vek{W}\vek{v} = 
    \diag(\vek{w})\vek{v} = \tns{M}(\tnb{w}) \tns{V}(\tnb{v}),
\end{equation}
which helps in understanding the following algorithms.
Note that this shows that the algebra $\bigotimes_{\ell=1}^d \D{R}^{M_\ell}$ 
with the Hadamard product is already jointly diagonalised
and is essentially in its \emph{Gel'fand} representation \citep{segalKunze78},
which is an abstract way of saying what was already stated above, namely that
the spectrum $\sigma(\tnb{w})$ of an element $\tnb{w}\in \bigotimes_{\ell=1}^d \D{R}^{M_\ell}$
are exactly the individual
terms or data stored in the tensor.  If in some abstract unital associative
and commutative algebra the multiplication is \emph{not} the point-wise
multiplication --- e.g.\ think of convolution --- then it is known that
modulo some technicalities  \citep{segalKunze78} it is isomorphic via the
Gel'fand ``diagonalisation'' morphism to a function algebra,
and all the algorithms to follow would deal with the spectrum 
$\sigma(\tnb{w})$ of $\tnb{w}$ instead of with the ``values of $\tnb{w}$''.
These two notions only coincide for a function algebra.

\subsection{Post-processing Tasks} \label{SS:post-proc-tasks}
Here the different post-processing tasks are explained together
with how they will be computed.  This may involve a number
of auxiliary functions.  The computation of these --- again through
truncated iteration  Algorithm~\ref{alg:basic} ---  and the way
in which some of the computations can be enhanced or accelerated
is shown in \refSS{auxil}.

\paragraph{Finding the maximum or minimum} of 
$\tnb{w}\in\C{T}=\bigotimes_{\ell=1}^d \D{R}^{M_\ell}$
is the first task we consider.  Observe, that the element of maximum
modulus of the tensor $\tnb{w}$ is also equal to the $\infty$-norm $\nd{\tnb{w}}_\infty$.

Finding the maximum means
finding the index $\vhat{m} = (\hat{m}_{1},\dots,\hat{m}_{d})\in\C{M}$
where the maximum $\hat{\tns{w}}$ of the elements in $\tnb{w}$ occurs:
\begin{equation}  \label{eq:max-loc}
   \hat{\tns{w}}:=\tns{w}_{\vhat{m}} := \tns{w}_{\hat{m}_{1},\dots,\hat{m}_{d}} := 
       \max \, \{\, \tns{w}_{\vek{m}} \; : \; \vek{m}=(m_1,\dots,m_d)\in\C{M}\} .
\end{equation}

It was already established that each value of $\tnb{w}$ is an eigenvalue.
Defining for each $1\le m\le M_\ell$ the canonical unit basis vectors 
$\vek{e}^{(m)}_{M_\ell}$ in each $\D{R}^{M_\ell}$ in the tensor product 
$\C{T}=\bigotimes_{\ell=1}^d \D{R}^{M_\ell}$ as 
$\vek{e}^{(m)}_{M_\ell}:= (\updelta_n^{(m)})_{n=1,\dots,M_\ell} \in \D{R}^{M_\ell}$ via 
the Kronecker-$\updelta$-symbol, and similarly for each $\vek{m}\in\C{M}$ the canonical unit 
basis vectors $\tnb{e}^{(\vek{m})}=(\updelta^{(\vek{m})}_{(\vek{n})})_{\vek{n}\in\C{M}} 
\in \C{T}=\bigotimes_{\ell=1}^d \D{R}^{M_\ell}$,
then for \emph{any} $\vek{m}=(m_1,\dots,m_d)\in\C{M}$ the element $\tns{w}_{\vek{m}}$
satisfies an \emph{eigenvalue} equation: 
\begin{equation} \label{eq:ev-max}
\tnb{L}_{\tnb{w}}
   \tnb{e}^{(\vek{m})} := \tns{V}^{-1}\left(\tns{M}(\tnb{w})\tns{V}(\tnb{e}^{(\vek{m})})\right)=
    \tnb{w} \odot \tnb{e}^{(\vek{m})} =
  \tns{w}_{\vek{m}} \tnb{e}^{(\vek{m})},
\end{equation}
and the eigenvectors $\tnb{e}^{(\vek{m})} = \bigotimes_{\ell=1}^d \vek{e}^{(m_\ell)}_{M_\ell}$
are evidently of rank one.
This shows that each datum $\tns{w}_{\vek{m}}$ of $\tnb{w}$ may be found through 
\emph{eigenvalue} computation.  The eigenvalue is the data entry
$\tns{w}_{\vek{m}}$ with the corresponding eigenvector $\tnb{e}^{(\vek{m})}$
indicating the location resp.\ the index $\vek{m}\in\C{M}$.

Although quite obvious, it may be worthwhile pointing out that the linear span of any
nonempty collection of eigen- resp.\ unit vectors
$\C{T}_{\C{N}} := \spn \{ \tnb{e}^{(\vek{n})}\mid \vek{n}\in\C{N} \}$
($\C{N}\subseteq\C{M}$) is an invariant subspace, and the orthogonal
projector on $\C{T}_{\C{N}}$ is $\tnb{L}_{\tnb{p}^{(\C{N})}} := \sum_{\vek{n}\in\C{N}}
 \tnb{L}_{\tnb{e}^{(\vek{n})}}$ with 
$\tnb{p}^{(\C{N})} = \sum_{\vek{n}\in\C{N}} \tnb{e}^{(\vek{n})}$.
With this, the spectral resolution of the identity for any and all 
$\tnb{L}_{\tnb{w}}\in\E{L}(\C{T})$
for $\tnb{w}\in\C{T}$ is $\tnb{I}_{\C{T}} = \sum_{\vek{m}\in\C{M}} 
\tnb{L}_{\tnb{e}^{(\vek{m})}}$ corresponding to
$\tnb{1} = \sum_{\vek{m}\in\C{M}} \tnb{e}^{(\vek{m})}$, giving the spectral resolution
\begin{equation} \label{eq:spec-dec}
\tnb{L}_{\tnb{w}} = \sum_{\vek{m}\in\C{M}} \tns{w}_{\vek{m}} \tnb{L}_{\tnb{e}^{(\vek{m})}},
\; \text{ corresponding to } \;
\tnb{w} = \sum_{\vek{m}\in\C{M}} \tns{w}_{\vek{m}} \tnb{e}^{(\vek{m})} = 
  \sum_{\vek{m}\in\C{M}} \tns{w}_{\vek{m}} \bigotimes_{\ell=1}^d \vek{e}^{(m_\ell)}_{M_\ell},
\end{equation}
which may be seen as the trivial basis representation of $\tnb{w}$.

As the basic algorithms for eigenvalue computation like power iteration
for self-adjoint linear maps converge to the eigenvalue of maximum absolute
value or modulus, it is now clear that if we assume that the maximum is also
the element of maximum absolute value, we only have to perform power
iteration.  This gives us the value $\hat{\tns{w}}=\tns{w}_{\vhat{m}}$,
and the location through its corresponding eigenvector $\tnb{e}^{(\vhat{m})}$.
How to find the minimum in this case, or what to do when the element
of maximum modulus is the minimum will be explained later in \refSS{auxil};
it will all be accomplished with spectral shifts familiar from eigenvalue
computations.

\paragraph{Finding the maximum or minimum of some function $f$} of $\tnb{w}$,
where it is assumed that $f: \sigma(\tnb{w})\to\D{R}$,
is not very difficult if one can compute $f(\tnb{w})$, which will be simply
understood as $f(\tnb{w}):= (f(\tns{w}_{\vek{m}}))_{\vek{m}\in\C{M}}$.
This definition of $f(\tnb{w})$ is also the one which comes via the
isomorphy with the algebra of diagonal matrices
$\F{diag}(\D{R},N) \subset \F{gl}(\D{R},N)$, where $f(\tnb{w})$ is
defined \citep{NHigham} through the corresponding matrix function
$f(\tns{M}(\tnb{w})) = f(\vek{W})$, and is also the same as the
one which comes from the abstract functional calculus on the Banach
algebra $\C{T}$ \citep{segalKunze78}.

Observe that one may compute $f(\tnb{w})$ for any real or complex function $f$ defined
on the spectrum $\sigma(\tnb{w})$ --- which is a finite set --- at least in principle,
through an interpolating polynomial $p(w)=\sum_k \alpha_k w^k$.   As usual, the number
$w\in\sigma(\tnb{w})$ is replaced by $\tnb{w}$ (or $\tnb{L}_{\tnb{w}}$) in the polynomial,
resulting in $f(\tnb{w}) = p(\tnb{w}) = \sum_k \alpha_k \tnb{w}^{\odot k}$.  This may
be highly inefficient, as in general the degree of the polynomial would be equal
to the number of data, in this case the huge number $N-1$.  For the functions which
will be needed for the post-processing tasks, the computation of the function
will be done differently, namely iteratively via Algorithm~\ref{alg:basic}.

Assume now that $f(\tnb{w})$ has been computed,
at least approximately.  To find maxima or minima of $f(\tnb{w})$,
one now has to simply apply the considerations of the preceding paragraph
to the tensor $f(\tnb{w})$.
Obviously, $f(\tnb{w})$ has the same spectral resolution as $\tnb{w}$ in \refeq{eq:spec-dec}:
\begin{multline} \label{eq:spec-dec-f}
f(\tnb{L}_{\tnb{w}}) = \sum_{\vek{m}\in\C{M}} f(\tns{w}_{\vek{m}}) \tnb{L}_{\tnb{e}^{(\vek{m})}},
\; \text{ corresponding to } \\
f(\tnb{w}) = \sum_{\vek{m}\in\C{M}} f(\tns{w}_{\vek{m}}) \tnb{e}^{(\vek{m})} = 
  \sum_{\vek{m}\in\C{M}} f(\tns{w}_{\vek{m}}) \bigotimes_{\ell=1}^d \vek{e}^{(m_\ell)}_{M_\ell}.
\end{multline}

\paragraph{Finding the index and value closest to a given number} $\rho\in\D{R}$
is now simply finding the eigenvector and eigenvalue of 
\begin{equation} \label{eq:eigv-close}
(\tnb{L}_{\tnb{w}} - \rho\, \tnb{I}_{\C{T}})^{-1} = \tnb{L}_{\tnb{y}} \quad
  \text{with} \quad \tnb{y} = (\tnb{w} - \rho\, \tnb{1})^{\odot-1} .
\end{equation}
This is a special case of the preceding paragraph for the function
$f:t \mapsto (t-\rho)^{-1}$.  Therefore,
if $\lambda_{\tnb{w}}$ is any eigenvalue of $\tnb{L}_{\tnb{w}}$, then the corresponding
eigenvalue of $\tnb{L}_{\tnb{y}}$ is the transformed one 
$\lambda_{\tnb{y}} = (\lambda_{\tnb{w}} - \rho)^{-1}$.  Hence the value of $\tnb{w}$
closest to $\rho$ is the eigenvalue of maximum modulus of $\tnb{L}_{\tnb{y}}$, and
the element of largest magnitude of $(\tnb{w} - \rho\, \tnb{1})^{\odot-1}$.
To find the element of smallest magnitude, or even the vanishing elements
of $\tnb{w}$, one would use $\rho=0$.
Thus one may use the same algorithms as in the maximum search above.  Observe that this
operation here requires the Hadamard inverse, which will be one of the auxiliary functions.

\paragraph{Finding the indices in a level set} requires
the so-called \emph{sign} function, where
$\sign: \C{T} \to \C{T}$ is defined component-wise as 
\begin{equation} \label{eq:sign-def}
\C{T}\ni(\sign(\tnb{w})_{m_1,\dots,m_d}) :=
      \begin{cases}
          \phantom{-}1, & \text{ if } \tns{w}_{m_1,\dots,m_d}  > 0; \\
          \phantom{-}0, & \text{ if } \tns{w}_{m_1,\dots,m_d}  = 0; \\
          -1, & \text{ if } \tns{w}_{m_1,\dots,m_d}  < 0 .
      \end{cases} 
\end{equation}
Again this is a special case of the previous one with a general function $f$.

The \emph{characteristic} function for a subset $S\subset\D{R}$
--- we shall only look at intervals --- $\chi_S: \C{T} \to \C{T}$,
is again defined component-wise:
\begin{equation} \label{eq:characteristic-def}
\C{T}\ni(\chi_S(\tnb{w})_{m_1,\dots,m_d}) :=
      \begin{cases}
          1, & \text{ if } \tns{w}_{m_1,\dots,m_d}  \in S; \\
          0, & \text{ if } \tns{w}_{m_1,\dots,m_d}  \notin S .
      \end{cases} 
\end{equation}
With these two auxiliary functions, it is possible to define the characteristic
function of a level set, i.e.\ all values between $\omega_1, \omega_2 \in \D{R}$.
We then have with $-\infty < \omega_1 < \omega_2 < \infty$:
\begin{equation} \label{eq:char-def}
   (\chi_S(\tnb{w}))_{m_1,\dots,m_d}) :=
      \begin{cases}
          \frac{1}{2}\,(\tnb{1} + \sign(\omega_2\,\tnb{1} - \tnb{w})), 
                  & \text{ if } S = ] -\infty, \omega_2[; \\
          \frac{1}{2}\,(\tnb{1} - \sign(\omega_1\,\tnb{1} - \tnb{w})),
                  & \text{ if } S = ]\omega_1, +\infty[; \\
          \frac{1}{2}\,(\sign(\omega_2\,\tnb{1} - \tnb{w}) - \sign(\omega_1\,\tnb{1}  -\tnb{w})),
                  & \text{ if } S = ]\omega_1, \omega_2[; 
      \end{cases} 
\end{equation}
Each case is easily computed with the sign function from \refeq{eq:sign-def}.
Hence, the indices of all components 
$\tensor*{\tns{w}}{_{m_1}_{\dots}_{m_d}} \in ]\omega_1\,\omega_2[$ are provided by
$\chi_{]\omega_1, \omega_2[}(\tnb{w})$.
One may additionally define the \emph{level set} function $\C{L}_S(\cdot)$ of a 
subset $S\subset \D{R}$ as $\C{L}_S(\tnb{w}) := \chi_S(\tnb{w}) \odot \tnb{w}$,
which together with the indices provides the values in the subset $S$.

\paragraph{Finding the number of indices in a level set} is accomplished
through consideration of the \emph{support} 
\begin{equation}  \label{eq:def-supp}
\text{supp}\,\chi_S := \{\vek{m}\in\C{M}\;:\; \chi_S(\tnb{w})_{\vek{m}} \neq 0 \} 
   \subset\C{M}
\end{equation}
of a characteristic function $\chi_S$ as the subset of those indices where it
is \emph{non-zero}.  Its \emph{cardinality} is 
\begin{equation}  \label{eq:def-card}
\#(\text{supp}\,\chi_S(\tnb{w})) = 
  |\text{supp}\,\chi_S(\mat{w})| = \bkt{\chi_S(\tnb{w})}{\tnb{1}}_{\C{T}},
\end{equation}
the number of non-zero positions in the characteristic function;
requiring only the computation of one inner product with the Hadamard multiplicative unit.

\paragraph{Computing the probability of being in a level, as well
as the mean and the variance} is under
the assumption that each index in $\tnb{w}$ carries the same probability,
as then the \emph{probability}
of a value of $\tnb{w}$ being in a subset $S$ is simply
\begin{equation}  \label{eq:def-prob}
   \D{P}_{\tnb{w}}(S) := \frac{\#(\text{supp}\,\chi_S(\tnb{w}))}{N}
     = \frac{\bkt{\chi_S(\tnb{w})}{\tnb{1}}_{\C{T}}}{N}.
\end{equation}

The  \emph{mean} or \emph{average} and \emph{variance} of $\tnb{w}$ is then  
\begin{equation}  \label{eq:def-mean-var}
 \EXP{\tnb{w}}:=\bar{\tns{w}} := \frac{1}{N}\, \bkt{\tnb{w}}{\tnb{1}}_{\C{T}}; \quad
  \text{var}(\tnb{w}) = \frac{1}{N}\, \bkt{\tnb{\tilde{w}}}{\tnb{\tilde{w}}}_{\C{T}},
  \; \text{ where } \tnb{\tilde{w}} := \tnb{w} - \bar{\tns{w}}\, \tnb{1}.
\end{equation}
Given the level set function from above, one may also compute the conditional mean,
conditioned on being in the set $S$:
\begin{equation}  \label{eq:def-mean-cond}
   \EXP{\tnb{w}|S} := \bar{\tns{w}}_{| S} := 
   \frac{\bkt{\C{L}_S(\tnb{w})}{\tnb{1}}_{\C{T}}}{\bkt{\chi_S(\tnb{w})}{\tnb{1}}_{\C{T}}} .
\end{equation}

\subsection{Auxiliary functions and algorithmic details}   \label{SS:auxil}
Here we show how to compute the auxiliary functions used in
the previous \refSS{post-proc-tasks}.  These are functions defined
on the Euclidean Hadamard algebra  $\C{T}=\bigotimes_{\ell=1}^d \D{R}^{M_\ell}$,
and it is worth while pointing out that the iterative algorithms
will be operating on the algebra, which, as we saw
in \refSS{prelim}, is isomorphic to an algebra of diagonal matrices,
$\F{diag}(\D{R},N) \subset \F{gl}(\D{R},N)$, and these matrices are real symmetric.
This means in particular that all the known algorithms for computing
matrix functions \citep{NHigham} can be used on the Hadamard algebra as well.
Of those needed here, the most basic one and the one with the
simplest connection to the algebra turns out to be the inverse $\tnb{w}^{\odot-1}$.
Let us remark once again that the algorithms are valid in any abstract algebra,
but we are mainly concerned with the example of the ``function algebra''
$\bigotimes_{\ell=1}^d \D{R}^{M_\ell}$.

\paragraph{The Hadamard inverse} $\tnb{w}^{\odot-1}$ of a Hadamard invertible
$\tnb{w}\in\C{T}$ is needed for \emph{inverse} iteration and other post-processing tasks.
Although the Hadamard inverse is the application of the function $f:t\mapsto t^{-1}$
to $\tnb{w}$ and could in principle be computed through a polynomial,
it is often advantageous and much more effective to use other algorithms.
The algorithm for the inverse can be given simply by referring to the
quadratically convergent Newton algorithm
for matrices  \citep{NHigham} for computing inverses, i.e.\ we
use \emph{Newton's method} for ($P\leftarrow (\odot -1)$) and
apply it to the equation $\tns{F}(\tnb{v}):=\tnb{v}^{\odot -1} - \tnb{w} = \tnb{0}$
to solve for $\tnb{v}$, the only solution of which is $\tnb{v}=\tnb{w}^{\odot -1}$.
Hence the iteration function to be used in Algorithm~\ref{alg:basic} with starting
vector $\tnb{v}_0 := \tnb{w}$ is
\begin{equation}  \label{eq:iter-Had-inv}
\tns{\Phi}_{(\odot -1)}(\tnb{v}) := \tnb{v}\odot(2 \cdot \tnb{1} - \tnb{w}\odot\tnb{v}).
\end{equation}
It is well known that the iteration converges \emph{quadratically} \citep{NHigham},
and according to what was explained in \refSS{itrunc} this is true
even with \emph{truncation}.

\paragraph{The sign function} is defined in \refeq{eq:sign-def}.
Many of the tasks explained in the previous \refSS{post-proc-tasks}
involve the sign function.    Given the sign
function, the characteristic, support, and level set functions are
easily computed by the basic operations of the algebra.

From the fact that each value in a tensor $\tnb{w}$ is an eigenvalue
in the Hadamard algebra, and from the isomorphy with the algebra of
diagonal matrices $\F{diag}(\D{R},N) \subset \F{gl}(\D{R},N)$
where each value on the diagonal is also obviously an eigenvalue,
one sees that this is actually the same definition as the one for
matrices \citep{NHigham}, and is also the one which follows from general
functional calculus in abstract algebras \citep{segalKunze78}.
This means that $\sign(\tnb{w})$ is the application of the function
$f:t\mapsto \sign(t)$ to $\tnb{w}\in\C{T}$, and again could in principle
be computed with a polynomial.
To compute it via the iterative Algorithm~\ref{alg:basic},
i.e.\ (\emph{$P\leftarrow\sign$}), one uses the same algorithm \citep{NHigham} ---
the Roberts-Newton algorithm --- as for matrices.
It can be derived by applying \emph{Newton's} method to the equation
$\tns{F}(\tnb{v}):=\tnb{v}\odot\tnb{v} - \tnb{1}=\tnb{0}$ with
starting value $\tnb{v}_0 := \tnb{w}$.
This yields the iteration function 
\begin{equation}  \label{eq:iter-sign-f-p}
\tns{\Phi}_{\sign}(\tnb{v}) := \frac{1}{2}\,(\tnb{v}+ \tnb{v}^{\odot -1}).
\end{equation}
Observe that this iteration function in a slightly more general form
$\tns{\Phi}_{(\sqrt{\tnb{w}})}(\tnb{v}) := 
\frac{1}{2}\,(\tnb{v}+ \tnb{v}^{\odot -1}\odot\tnb{w})$ is the ancient
\emph{Babylonian} method  \citep{NHigham} --- in modern paralance Newton's method --- to
find the square root $\tnb{w}^{1/2}$ of $\tnb{w}$, as by
inserting $\tnb{w}=\tnb{1}$ in the Babylonian iteration $\tns{\Phi}_{(\sqrt{\tnb{w}})}$,
one obtains \refeq{eq:iter-sign-f-p}.  This means that one is iterating
to compute the square root of the unit element $\tnb{1}$ with the specific
starting value $\tnb{v}_0 = \tnb{w}$.  We will
come back to this point of view shortly.

It is known \citep{NHigham} that the iteration converges \emph{quadratically},
and according to what was explained in \refSS{itrunc}, 
it does so even with \emph{truncation}; observe that it needs
the Hadamard inverse from the previous paragraph.  Therefore it is 
not really practical as this would mean that it can only be applied
to invertible $\tnb{w}\in\C{T}=\bigotimes_{\ell=1}^d \D{R}^{M_\ell}$, but it also
actually needs a new inverse $\tnb{v}_i^{\odot -1}$ in each iteration.  This last 
fact would computationally result in a nested iteration --- using 
the algorithm from the previous paragraph in an inner iteration to
compute  the inverse $\tnb{v}_i^{\odot -1}$ in each sweep of the
outer iteration for the sign function ---
and such procedures are seldom computationally efficient.

It is thus simpler to use a device employed also for matrix sign
computations  \citep{NHigham}, namely to replace the explicit
inverse in \refeq{eq:iter-sign-f-p} by one step of the iteration
for the inverse in \refeq{eq:iter-Had-inv}.  The resulting
Newton-Schulz algorithm \citep{NHigham} has the same starting
point $\tnb{v}_0 := \tnb{w}$ as before, but uses the iteration function
\begin{equation}  \label{eq:iter-sign-N-S}
\tns{\Phi}_{\text{N-S}}(\tnb{v}) := 
   \frac{1}{2}\,\cdot \tnb{v}\odot(3\,\cdot \tnb{1} - \tnb{v}^{\odot 2})
\end{equation}
in the Algorithm~\ref{alg:basic}.  It is now a simple iteration which still converges
\emph{quadratically} \citep{NHigham}, even with \emph{truncation}
according to the explanations in \refSS{itrunc}.
We refer to the discussion in \citep{NHigham} on how to get even faster
algorithms using Padé approximations and other ways to accelerate
the convergence through scaling, which is especially important in
the initial stages.

Coming back to the somewhat curious idea of iterating for the square
root of $\tnb{1}$, one may observe that although
the defining equation $\tns{F}(\tnb{v}):=\tnb{v}\odot\tnb{v} - \tnb{1}=\tnb{0}$
does not cover the case that a datum with the exact value of zero occurs ---
a non-invertible element would not satisfy the defining equation --- the Newton-Schulz
iteration with the function in \refeq{eq:iter-sign-N-S} with starting value
$\tnb{v}_0 := \tnb{w}$ takes care of this.  A value which vanishes
in $\tnb{v}_i$ also vanishes in $\tnb{v}_{i+1}=\tns{\Phi}_{\text{N-S}}(\tnb{v}_i)$,
due to the product with $\tnb{v}_i$ in \refeq{eq:iter-sign-N-S}.

\paragraph{Stopping criteria in case of quadratic convergence} 
in an iteration like in Algorithm~\ref{alg:basic}
are well known.  First, if one wants to solve $\tns{F}(\tnb{v})=\tnb{0}$
or a fixed point equation $\tns{F}(\tnb{v}):=\tns{\Phi}(\tnb{v})-\tnb{v}=\tnb{0}$,
a natural criterion is the size of the residuum at step $i$:
\begin{equation}  \label{eq:stop-F-all}
    \nd{\tns{F}(\tnb{v}_i)}_{\C{T}} < \eta_{\tns{F}}.
\end{equation}
But this checks only how well the equation is satisfied, and not directly
how accurate the iterate $\tnb{v}_i$ is.  Regarding this latter issue,
specifically for quadratic convergence,  a natural criterion \citep{NHigham}
to check for the accuracy of $\tnb{v}_i$ at step $i$ is
\begin{equation}  \label{eq:stop-v-n}
\delta_{i} :=  
   \frac{\nd{\tnb{v}_i - \tnb{v}_{i-1}}_{\C{T}}}{\nd{\tnb{v}_{i}}_{\C{T}}} < \eta_{\tns{v}}.
\end{equation}
Further, referring to the discussion specifically for the sign function
in \citep{NHigham}, there are arguments to check
\begin{equation}  \label{eq:stop-v-p}
\delta_{i} <  \nd{\tnb{v}_{i}}_{\C{T}}^p \,\eta_{\tns{v}}
\end{equation}
for the exponents $p=0,1,2$.  For $p=0$, this is the original general
criterion \refeq{eq:stop-v-n}, whereas the other values of $p$ take 
specific consideration of the sign function.

\paragraph{Eigenvalue computations} are involved in the first three
tasks described in \refSS{post-proc-tasks}.  One may use any algorithm
developed for large scale matrices \citep{golubVanloan, parlett98, saad92, ds-watk-07}
which only uses the action of the matrix on a vector in the computation
--- this is the Hadamard product in our case ---
and we shall here only explain the main idea and a few variations.

The simplest algorithm is power iteration ($P\leftarrow\text{pow-it}$),
and the iteration map $\tns{\Phi}_{\text{pow-it}}$ is given in
Algorithm~\ref{alg:power-it}, to be used in Algorithm~\ref{alg:basic}.  
Assume that the datum $\tns{w}_{\vhat{m}}>0$
of $\tnb{w}$ with maximum absolute value or maximum modulus is indeed a (positive)
maximum.  Assume also that this occurs at one unique index $\vhat{m}$,
i.e.\ the datum $\lambda_1:=\tns{w}_{\vhat{m}}$ is a simple eigenvalue,
and denote the next by absolute value smaller element / eigenvalue by $\lambda_2$.
Assume the further eigenvalues ordered by decreasing absolute value.

\begin{algorithm}
  \caption{One step power iteration $\tns{\Phi}_{\text{pow-it}}$}
       \label{alg:power-it}
   \begin{algorithmic}[1]     
    \State Input iterate $\tnb{v}_i$;
    \Comment{ assume $\tnb{v}_i$ of unit length. }
    \State $\tnb{u} \gets \tnb{w} \odot \tnb{v}_i$;
    \Comment{ $\tnb{u} = \tnb{L}_{\tnb{w}} \tnb{v}_i$. }
    \State $\gamma \gets \bkt{\tnb{u}}{\tnb{u}}_{\C{T}}^{-1/2}$; 
    \Comment{ inverse length of $\tnb{u}.\quad\lambda_1 \approx \gamma^{-1}$ }
    \State $\tnb{z} \gets \gamma\;\cdot \tnb{u}$;
    \Comment{ normalise output $\tnb{z}$ to unit length. }
    \State Output $\tnb{z}$;
  \end{algorithmic}
\end{algorithm}

As is well known  \citep{golubVanloan, parlett98},
the generated sequence $\tnb{v}_i$ \emph{converges linearly} to 
$\pm\tnb{e}^{(\vhat{m})}$ with a contraction factor of $q = |\lambda_2 / \lambda_1|$,
and with what was said in the explanations in \refSS{itrunc},
it still converges with truncation and will stagnate in the
vicinity of $\pm\tnb{e}^{(\vhat{m})}$.
Also, as in the algorithm the input
$\tnb{v}_i$ has unit length, the length of $\tnb{u}$
--- given by the inverse $\gamma^{-1}$ of the scaling factor in line 3 ---
converges linearly to the desired $\lambda_1=\tns{w}_{\vhat{m}}$
with the same rate.  Algorithmically, the eigenvalue approximation
$\lambda_1 \approx \gamma^{-1}$ is a \emph{side-effect} of Algorithm~\ref{alg:power-it}.

This can be immediately enhanced through 
the additional computation of the \emph{Rayleigh quotient} (RQ)
\begin{equation}   \label{eq:RQ-def}
\vrho_{\tnb{w}}(\tnb{v}_i) := 
   \frac{\bkt{\tnb{L}_{\tnb{w}}\tnb{v}_i}{\tnb{v}_i}_{\C{T}}}{\bkt{\tnb{v}_i}{\tnb{v}_i}_{\C{T}}}
 = \frac{\bkt{\tnb{w}\odot\tnb{v}_i}{\tnb{v}_i}_{\C{T}}}{\bkt{\tnb{v}_i}{\tnb{v}_i}_{\C{T}}}
\end{equation}
between line 2 and line 3 of Algorithm~\ref{alg:power-it}: 
$\vrho_{\tnb{w}}(\tnb{v}_i) = \bkt{\tnb{u}}{\tnb{v}_i}_{\C{T}}$
--- no need to divide by the length of the unit vector $\tnb{v}_i$.
The Rayleigh quotient is stationary at an eigenvalue and has
thus a quadratic convergence, it usually
represents a much better approximation to $\lambda_1=\tns{w}_{\vhat{m}}$
than $\gamma^{-1}$; see Algorithm~\ref{alg:power-RQ}.

\paragraph{Stopping criteria for eigenvalues} could be done similarly
to normal convergent processes, but is is possible to have a posteriori
error estimates which are actual bounds specifically for eigenvalues
and -vectors, e.g.\ see \citep{hgm-eigenbounds-85, golubVanloan, parlett98}.
The simplest seems to be the so-called Krylov-Bogolyubov bound.
With an approximate eigenvalue $\mu$ and approximate eigenvector $\tnb{x}$
of $\tnb{L}_{\tnb{w}}$ it is given by \citep{hgm-eigenbounds-85}
\begin{equation} \label{eq:kryl-bogo}
   \min_{\lambda_j\in\sigma(\tnb{w})} \frac{\ns{\lambda_j-\mu}}{\ns{\lambda_j}}
   \le \frac{\nd{\tnb{L}_{\tnb{w}}\tnb{x} - \mu\,\tnb{x}}_{\C{T}}}{\nd{\tnb{x}}_{\C{T}}}.
\end{equation}
The right-hand side of \refeq{eq:kryl-bogo} is minimised by $\mu = \vrho_{\tnb{w}}(\tnb{x})$.  
A short computation \citep{hgm-eigenbounds-85} shows that for the substitutions
$\tnb{x}\gets\tnb{v}_i$ and $\mu \gets \vrho_{\tnb{w}}(\tnb{v}_i)$
the right-hand side of \refeq{eq:kryl-bogo} becomes $(\vrho_{\tnb{w}^{\odot 2}}(\tnb{v}_i)
- \vrho_{\tnb{w}}(\tnb{v}_i)^2)^{1/2}$.  As 
\[
\vrho_{\tnb{w}^{\odot 2}}(\tnb{v}_i) =
\bkt{\tnb{w}^{\odot 2}\odot\tnb{v}_i}{\tnb{v}_i}_{\C{T}} = 
\bkt{\tnb{w}\odot\tnb{v}_i}{\tnb{w}\odot\tnb{v}_i}_{\C{T}} = \bkt{\tnb{u}}{\tnb{u}}_{\C{T}},
\]
the Krylov-Bogolyubov error bound \refeq{eq:kryl-bogo} can be computed as
\begin{equation} \label{eq:kryl-bogo-2}
   \min_{\lambda_j\in\sigma(\tnb{w})} 
     \frac{\ns{\lambda_j-\bkt{\tnb{u}}{\tnb{v}_i}_{\C{T}}}}{\ns{\lambda_j}}
   \le \vepsilon_\lambda :=
   (\bkt{\tnb{u}}{\tnb{u}}_{\C{T}} - \bkt{\tnb{u}}{\tnb{v}_i}_{\C{T}}^2)^{1/2}.
\end{equation}
Observe that this is not an a posteriori error estimate, but an actual \emph{bound}.
It is the best possible with the available information \citep{hgm-eigenbounds-85}, and
it can be very easily computed inside the iteration to control the possible termination.
Similar bounds exist for the approximate eigenvector $\tnb{v}_i$.  Inserting these
considerations into Algorithm~\ref{alg:power-it} gives Algorithm~\ref{alg:power-RQ}
for the task power iteration with RQ computation ($P \gets \text{pow-RQ}$)
for the iteration map $\tns{\Phi}_{\text{pow-RQ}}$ to be used in Algorithm~\ref{alg:basic},
which has the more accurate RQ eigenvalue approximation $\lambda_1 \approx \vrho_{\tnb{w}}$
and error bound $\vepsilon_\lambda =
   (\bkt{\tnb{u}}{\tnb{u}}_{\C{T}} - \bkt{\tnb{u}}{\tnb{v}_i}_{\C{T}}^2)^{1/2}$
included as side-effects.

\begin{algorithm}
  \caption{One step power iteration with Rayleigh quotient (RQ): $\tns{\Phi}_{\text{pow-RQ}}$}
       \label{alg:power-RQ}
   \begin{algorithmic}[1]     
    \State Input iterate $\tnb{v}_i$;
    \Comment{ assume $\tnb{v}_i$ of unit length. }
    \State $\tnb{u} \gets \tnb{w} \odot \tnb{v}_i$;
    \Comment{ $\tnb{u} = \tnb{L}_{\tnb{w}} \tnb{v}_i$. }
    \State $\vrho_1 \gets \bkt{\tnb{u}}{\tnb{v}_i}_{\C{T}}$; 
    \Comment{ Rayleigh quotient (RQ) $\lambda_1 \approx \vrho_1$ }
    \State $\vrho_2 \gets \bkt{\tnb{u}}{\tnb{u}}_{\C{T}}$; 
    \Comment{ RQ of $\tnb{w}^{\odot 2}$. }
    \State $\vepsilon_\lambda \gets (\vrho_2-\vrho_1^2)^{-1/2}$; 
    \Comment{ error bound \refeq{eq:kryl-bogo-2} of $\lambda_1$. }
    \State $\gamma \gets \vrho_2^{-1/2}$; 
    \Comment{ inverse length of $\tnb{u}$. }
    \State $\tnb{z} \gets \gamma\;\cdot \tnb{u}$;
    \Comment{ normalise output $\tnb{z}$ to unit length. }
    \State Output $\tnb{z}$;
  \end{algorithmic}
\end{algorithm}

Once the determination
of the eigenvalue is accurate enough such that one has an interval which contains only
\emph{one} eigenvalue, one may use the even tighter
Temple-Kato bounds \citep{hgm-eigenbounds-85}.

\paragraph{Starting vectors and other enhancements} like deflation and
Krylov subspaces are discussed here.
Recalling the discussion of invariant subspaces in the paragraph on maxima
and minima in \refSS{post-proc-tasks}, the starting vector $\tnb{v}_0$
in Algorithm~\ref{alg:basic} should not have
any zero in it, as this would exclude an invariant subspace from the
investigation.  One possibility is $\tnb{v}_0 = \tnb{1}/\sqrt{N}$, being equal
in all positions.  With this starting vector one has $\tnb{v}_1 \propto 
\tnb{w}\odot\tnb{1}=\tnb{w}^{\odot 1}$, and it is easy to see that
what is computed are scaled versions of $\tnb{w}^{\odot i}$ 
as eigenvector approximation $\tnb{v}_i$.  This gives another possibility
for acceleration \citep{Grasedyck-gamm}, which will be discussed later.

Another technique in eigenvalue computations is \emph{deflation}
\citep{golubVanloan, parlett98}: when one eigenvalue $\lambda_{1}=\lambda_{\vhat{m}}$
of largest modulus and corresponding eigenvector $\tnb{e}^{(\vhat{m})}$
have been located, one may want to compute a further location $\vtil{m}$ where
there is an element $\tns{w}_{\vtil{m}}=\lambda_{\vtil{m}}$ of equal magnitude 
$|\tns{w}_{\vtil{m}}|=|\lambda_{1}|=|\lambda_{\vhat{m}}|$.  For this one can use
deflation, and one could change line 2 in Algorithm~\ref{alg:power-it} or 
Algorithm~\ref{alg:power-RQ} to $\tnb{u}\gets \tnb{w}\odot (\tnb{1} - \tnb{e}^{(\vhat{m})})
\odot \tnb{v}_i$ to obtain a new iteration map for deflation.
The factor $(\tnb{1} - \tnb{e}^{(\vhat{m})})$ projects into the
orthogonal complement of the invariant subspace $\spn\{ \tnb{e}^{(\vhat{m})}\}$,
and thus all the eigenvalues are left unchanged, except for $\lambda_{\vhat{m}}$,
which is mapped to zero.  But in this special case --- where we know the form
of all invariant subspaces --- it is even simpler to keep the same iteration
map and to choose
as starting vector $\tnb{v}_0:=(\tnb{1} - \tnb{e}^{(\vhat{m})})/\sqrt{(N-1)}$.  It has a
zero at position $\vhat{m}$ and is thus in the invariant subspace 
$(\spn\{\tnb{e}^{(\vhat{m})}\})^\perp$, and all iterates will stay in that subspace.
Iterating in Algorithm~\ref{alg:basic} with that starting vector either
with the iteration map in Algorithm~\ref{alg:power-it} or in
Algorithm~\ref{alg:power-RQ} would give us either another eigenvalue $\lambda_{\vtil{m}}$ 
of equal magnitude, or the eigenvalue with second largest magnitude $\lambda_2$,
and the corresponding eigenvector.  In this manner all desired eigenvalues can be
computed via deflation. 

Some other acceleration techniques should be mentioned briefly:  
The convergence speed in power iteration is
controlled by the ratio $|\lambda_2 / \lambda_1|$, where $\lambda_2$
is the next smallest eigenvalue in absolute size.  Sometimes this ratio can be very close
to unity.  If instead with a single vector $\tnb{v}_i$ one iterates with
a whole block $\tnb{v}_i^{(1)},\dots,\tnb{v}_i^{(j)}$ of $j$ mutually orthogonal
vectors --- effectively a subspace --- the convergence speed changes to
$|\lambda_{1+j} / \lambda_1|$.  One has to restore orthogonality though after
each iteration sweep \citep{parlett98} of this block- or subspace iteration.
This kind of technique also makes it possible to compute multiple eigenvalues,
i.e.\ when the maximum occurs at several places.

Other well-known methods for eigenvalues with even faster convergence
build on \emph{Krylov} subspaces
\citep{parlett98}, here we have symmetric matrices and thus one would
use the Lanczos method.  These procedures rely on orthogonalisation though,
and this may be problematic when combined with truncation.  And certainly
the Krylov subspace can be combined with the block iteration idea to
give block- or subspace-Lanczos methods \citep{hgm-subLanczos85}.

\paragraph{``Exponentiating'' the power iteration} can be easily achieved
via the clever idea of \citep{Grasedyck-gamm}. Recall that with the starting
vector $\tnb{v}_0:=\tnb{1}/\sqrt{N}$, power iteration would compute 
$\tnb{v}_1 \propto \tnb{w}^{\odot 1}=\tnb{w}^{\odot 1}\odot \tnb{1}$,
a scaled version of $\tnb{w}$,
corresponding to the action of $\tnb{w}^{\odot 1}$ on $\tnb{1}$.  
Starting actually with $\tnb{v}_0\propto\tnb{w}$ and 
changing line 2 in Algorithm~\ref{alg:power-it} or Algorithm~\ref{alg:power-RQ}
to $\tnb{u} := \tnb{v}_i^{\odot 2}$, one has $\tnb{v}_1 \propto \tnb{w}^{\odot 2}=
\tnb{w}^{\odot 2}\odot\tnb{1}$, i.e.\ the action of $\tnb{w}^{\odot 2}$ on $\tnb{1}$.
In the next iteration
one has $\tnb{v}_2 \propto \tnb{w}^{\odot 4}$, and in iteration $i$
one has $\tnb{v}_i \propto \tnb{w}^{\odot 2^i}$.
Reformulating Algorithm~\ref{alg:power-RQ} for this new task of
``exponentiated'' power iteration ($P\gets\text{exp-pow}$)
for the new iteration function $\tns{\Phi}_{\text{exp-pow}}$ results 
in Algorithm~\ref{alg:exp-pow}, to be used in Algorithm~\ref{alg:basic}
with starting vector $\tnb{v}_0 := \tnb{w}/\nd{\tnb{w}}_{\C{T}}$.
We also need to keep the auxiliary vector $\tnb{y}=\tnb{w}^{\odot 2}$.

\begin{algorithm}
  \caption{One step ``exponentiated'' power iteration with RQ: $\tns{\Phi}_{\text{exp-pow}}$}
       \label{alg:exp-pow}
   \begin{algorithmic}[1]     
    \State Input iterate $\tnb{v}_i$;
    \Comment{ assume $\tnb{v}_i \propto \tnb{w}^{\odot 2^{i}}$ of unit length. }
    \State $\tnb{u} \gets \tnb{v}_i \odot \tnb{v}_i$;
    \Comment{ $\tnb{u} \propto \tnb{L}_{\tnb{w}}^{\odot 2^{i+1}} \tnb{1}=\tnb{w}^{\odot 2^{i+1}}$.}
    \State $\vrho_1 \gets \bkt{\tnb{w}}{\tnb{u}}_{\C{T}}$; 
    \Comment{ Rayleigh quotient (RQ) $\lambda_1 \approx \vrho_1$ }
    \State $\vrho_2 \gets \bkt{\tnb{y}}{\tnb{u}}_{\C{T}}$; 
    \Comment{ RQ of $\tnb{w}^{\odot 2}$. }
    \State $\vepsilon_\lambda \gets (\vrho_2-\vrho_1^2)^{-1/2}$; 
    \Comment{ error bound \refeq{eq:kryl-bogo-2} of $\lambda_1$. }
    \State $\gamma \gets \bkt{\tnb{u}}{\tnb{u}}_{\C{T}}^{-1/2}$; 
    \Comment{ inverse length of $\tnb{u}$. }
    \State $\tnb{z} \gets \gamma\;\cdot \tnb{u}$;
    \Comment{ normalise output $\tnb{z}$ to unit length. }
    \State Output $\tnb{z}$;
  \end{algorithmic}
\end{algorithm}

Thus the eigenvalue approximation is through the Rayleigh quotient
$\vrho_1 = \bkt{\tnb{w}}{\tnb{u}}_{\C{T}} = \bkt{\tnb{w}}{\tnb{v}_i^{\odot 2}}_{\C{T}}
= \bkt{\tnb{w}\odot\tnb{v}_i}{\tnb{v}_i}_{\C{T}}$ in line 3
of Algorithm~\ref{alg:exp-pow}.
As $\tnb{y}=\tnb{w}^{\odot 2}$, on line 4 the quantity
$\vrho_2= \bkt{\tnb{y}}{\tnb{u}}_{\C{T}} = 
\bkt{\tnb{w}^{\odot 2}}{\tnb{v}_i^{\odot 2}}_{\C{T}} = 
\bkt{\tnb{w}^{\odot 2}\odot\tnb{v}_i}{\tnb{v}_i}_{\C{T}}$ computes the RQ of 
$\tnb{w}^{\odot 2}$.  The ``eigenvector'' $\tnb{v}_i\propto\tnb{w}^{\odot 2^{i}}$
converges with the rate $|\lambda_2/\lambda_1|^{2^i}$, and hence one achieves
\emph{exponential} convergence.  This acceleration technique can be used in all
eigenvalue computations described here.

\paragraph{Transformed eigenvalue computations} are a well-known device
to find e.g.\ the maximum of the function $\tnb{y}=f(\tnb{w})$.  After
computing $f(\tnb{w})$, one proceeds as before, but with $\tnb{y}$
instead of $\tnb{w}$.  Some functions $f(t)$ are simple
enough so that the computation of $f(\tnb{w})$ is combined with the iteration.
Recall that this transforms all eigenvalues $\lambda_{\tnb{w}}$ to
$\lambda_{\tnb{y}}=f(\lambda_{\tnb{w}})$ according to spectral calculus,
leaving the eigenvectors unchanged.

\paragraph{Shifting and inverse shifting} are functions needed to access other than the
point of largest modulus in the spectrum of
$\tnb{L}_{\tnb{w}}$ resp.\ $\tnb{w}$.  These functions are the shift $f: t\mapsto t + \beta$
and the inverse shift $f:t \mapsto (t-\rho)^{-1}$.

The simple function of shifting $f(t) = t + \beta$ is needed
if the actual element of maximum modulus is a minimum of $\tnb{w}$;
then this is a maximum of $-\tnb{w}$.  In the iteration this will be picked up
by the Rayleigh-quotient being negative: $\vrho_{\tnb{w}} < 0$.
Also note that for this case --- where the maximum has smaller absolute
value than the minimum --- one may \emph{shift} the tensor by
a value of $\beta = - \vrho_{\tnb{w}} \approx -\lambda_1$ to
$\That{w}:=\tnb{w}+ \beta\cdot\tnb{1}$ so that the
maximum of $\That{w}$ has larger absolute value than the minimum of
$\That{w}$.  One then iterates with $\That{w}$.
After determining the maximum of $\That{w}$,
we may subtract the number $\beta>0$ again, to obtain the maximum of $\tnb{w}$.
The index of the maximum of $\That{w}$ is of course the same as the
one for the maximum of $\tnb{w}$.  Completely analogous manipulations can
be performed to find a minimum which is not of maximum modulus. 
These and similar techniques are well known from eigenvalue
calculations of large / sparse symmetric matrices  
\citep{golubVanloan, parlett98, saad92, ds-watk-07}.

The  inverse shift $f(t)=(t-\rho)^{-1}$ is needed to access intermediate points
in the spectrum $\sigma(\tnb{w})$.  Here one finds the data in $\tnb{w}$
closest to $\rho$, as they have the maximum modulus under the transformation.
Thus one computes first $\tnb{y}:=(\tnb{w} - \rho\cdot\tnb{1})^{\odot-1}$
approximately through Algorithm~\ref{alg:basic} with the iteration
function \refeq{eq:iter-Had-inv} --- with $\Tcek{w} := \tnb{w} - \rho\cdot\tnb{1}$
instead of $\tnb{w}$ --- followed by the eigenvector computation
with iteration function as in Algorithm~\ref{alg:exp-pow}, but with
starting vector $\tnb{v}_0 := \tnb{y}/\nd{\tnb{y}}_{\C{T}}$.

\newcommand{\CP}{\mrm{CP}}
\newcommand{\Tuck}{\mrm{T}}
\newcommand{\TT}{\mrm{TT}}

\section{Tensor formats}  \label{S:tensor-rep}
\subsection{Overview of tensor formats}
In this section we review definitions and properties of frequently used tensor formats.
Many such formats are used in quantum physics under the name \emph{tensor networks},
see \citep{VI03, Sachdev2010-a, EvenblyVidal2011-a, Orus2014-a, 
BridgemanChubb2017-a, BiamonteBergholm2017}.  We only look
at the \emph{canonical polyadic} (CP) and the \emph{Tensor Train}
(TT) format. For the Tucker \citep{Tuck:66} format see in \cite{Espig2019}.
The CP \citep{Hitch:27} format
have been well known for a long time and are therefore very
popular.  The TT format was originally developed in quantum physics and
chemistry as ``matrix product states'' (MPS), see \citep{VI03} and references therein, 
and rediscovered in \cite{oseledetsTyrt2010, oseledets2011} as tensor train.

A new class of tensor formats, which 
we do not consider, is the hierarchical 
tensor (HT) format.  It was introduced in \citep{HackbuschKuhn2009}, and further 
considered in \citep{Grasedyck2010}.

We note that the sets of low-rank tensors 
of fixed rank in the CP format are not closed for $d>2$, whereas
in the Tucker, TT, or HT tensor formats these sets are closed. The computations in the TT and HT formats
are based on the \emph{singular value decomposition} (SVD) \citep{OSTY09, Grasedyck2010}, 
and in the Tucker format on the  \emph{higher order SVD} (HOSVD) \citep{Tucker_trunk12}. 
In all tensor formats the tensor rank doubles for addition and squares for the
Hadamard product.  To avoid unnecessary and harmful rank growth,
the rank is usually truncated, e.g.\ by ALS-like or other optimisation algorithms 
\citep{EspigDiss, Espig2012VarCal}.

A tensor format is described by a
parameter vector space $\C{P}=\bigtimes_{\nu=1}^{d} \C{P}_\nu$,
where $\C{P}_\nu = \D{R}^{d_\nu}$, and a multilinear map 
$\tns{U}:\C{P} \rightarrow \C{T}$ into the tensor space 
$\C{T}:=\bigotimes_{\ell=1}^d \D{R}^{M_\ell}$.
For practical implementations of high dimensional problems we need
to distinguish between a tensor $\tnb{w} \in \C{T}$ and its tensor format
representation $\tnb{P} \in \C{P}$, where $\tnb{w}=\tns{U}(\tnb{P})$.  There are many
possibilities to define tensor formats.  Here, we consider the
canonical (CP), and the tensor train (TT) format.  In the following we
briefly repeat definitions and properties of these formats \citep{HA12}.
We also note that the reader can invent his own tensor format, which especially
well fits to his needs.

The CP format is cheap, it is simpler than the Tucker or TT format, but, 
compared to others, there are no reliable algorithms to compute CP decompositions 
for $d>2$ \citep{HA12, khorBook18}.
The Tucker format has stable algorithms \citep{Khoromskij_Low_Tacker}, 
but the storage and complexity costs are 
 $\C{O}(d\,r\,n +r^d)$, i.e.\ they grow exponentially with $d$. 
The TT format is a bit more complicated, but does not have this disadvantage.
CP and Tucker rank-structured tensor formats have been 
applied in chemometrics and in signal processing \citep{smilde-book-2004, Cichocki:2002}.

\subsection{The canonical polyadic tensor format}   \label{SS:CPFormat}
The canonical representation of multivariate functions  \citep{Hitch:27} was
introduced in 1927.
\begin{figure}[h]
\centering
\includegraphics[width=0.7\textwidth]{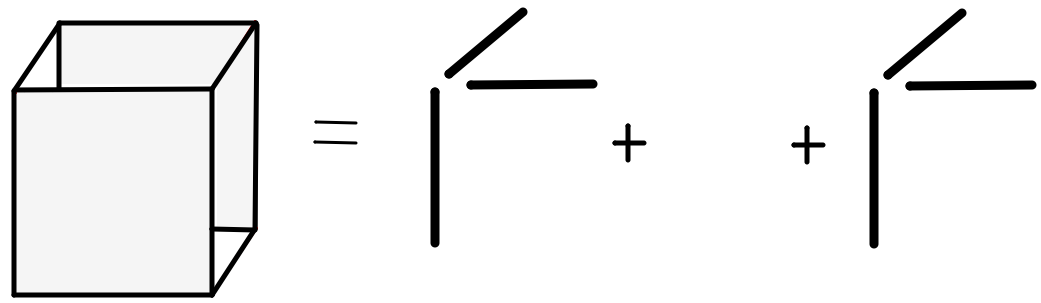}\quad\quad
\caption{Schema of the CP tensor decomposition of a 3D tensor.}
\label{fig:CP}
\end{figure}
A schema of the CP tensor format for $d=3$ is shown in 
Figure~\ref{fig:CP}. The full tensor 
$\tnb{w}\in \D{R}^{n_1\times n_2\times n_3}$ is shown on the left side of the equality 
sign, and its decomposition on the right.  The lines denote the vectors
$\tns{w}_{i 1}$, $\tns{w}_{i 2}$, and 
$\tns{w}_{i 3}$ for $i=1,\dots,r$, respectively.

\begin{definition}[Canonical Polyadic (CP) Tensor Format]\label{D:CPFormat}
The canonical polyadic tensor format in $\C{T}$ for variable $r$ is
defined by the multilinear mapping
\begin{eqnarray}    \label{eq:CP-format}
\tns{U}_{\CP,r} : \C{P}_{\CP,r} := \bigtimes_{\nu = 1}^d \C{P}^r_\nu &\rightarrow& \C{T},
     \quad \C{P}_\nu = \D{R}^{M_\nu},\\  \nonumber
     \C{P}_{\CP,r}  \ni
    \tnb{P}:= (\vek{w}_{i}^{(\nu)}: 1\leq i \leq r, \quad 1 \leq \nu \leq d) &\mapsto& 
    \tns{U}_{\CP,r}(\tnb{P}) := \tnb{w} = \sum_{i=1}^r \Ten_{\nu=1}^d \vek{w}_{i}^{(\nu)}\in\C{T}.
\end{eqnarray}
We call the sum of elementary tensors $\tnb{w}:= \tns{U}_{\CP,r}(\tnb{P})$
a tensor represented in the canonical tensor format with
$r$ terms.  The system of vectors $\tnb{P}= (\vek{w}_{i}^{(\nu)}\in\D{R}^{M_\nu} : 1\leq i \leq r,
\quad 1\leq \nu \leq d)$ is a \emph{representation system}
of $\tnb{w}$ with \emph{representation rank} $r$.  One may think of $\tnb{P}$
as of a vector valued $r\times d$ matrix with the vector $\vek{w}_{i}^{(\nu)}\in\D{R}^{M_\nu}$
at index position $(i,\nu)$.
\end{definition}
The storage requirement for $\tnb{w}=\tns{U}_{\CP,r}(\tnb{P})$ 
is $r\times \sum_{\nu=1}^d M_{\nu}$, and
in the simple case $M_1 = \dots = M_d = n$ it is $\C{O}(r\,d\,n)$.

\subsubsection{Basic operations with the canonical format}  \label{SSS:CP-ops}
%
%

We denote the set of all tensors $\tnb{w} \in \C{T}$ of rank $r$ by $\C{T}^r$.  The set
$\C{T}^r$ is a cone, i.e.\ $\tnb{w} \in \C{T}^r$ implies $\alpha\cdot \tnb{w} \in \C{T}^r$
for $\alpha\in\D{R}$, and at the same time $\C{T}^r$ is not a vector space as
for $\tnb{w}_1, \tnb{w}_2 \in  \C{T}^r $ one has in general 
$\tnb{w}_1+\tnb{w}_2 \notin  \C{T}^r $, but 
$\tnb{w}_1+\tnb{w}_2 \in  \C{T}^{2r} $ \citep{Espig2012VarCal}.


A complete description of fundamental operations in the canonical
tensor format and a their numerical cost can be found in
\citep{HA12}. For recent algorithms in the canonical tensor format we
refer to \citep{EspigDiss, ESHAGA09, ESHA09_1, ESHAROSCH09_1}.
%

\paragraph{Multiplication by a scalar}
$\alpha\in\D{R}$ could be done for a
$\tnb{w}=\tns{U}_{\CP,r}(\tnb{P})$ as in \refeq{eq:CP-format}
by multiplying all of the vectors $\{\vek{w}_{j}^{(\nu)},i=1,\dots,r\}$ 
for any $\nu$ by $\alpha$.  But to spread
the effect equally and keep the vectors balanced in size, we recommend
to define $\alpha_{\nu} := \sqrt[d]{|\alpha|}$ for all $\nu > 1$,
and $\alpha_{1} := \sign(\alpha)\sqrt[d]{|\alpha|}$.  Then, with a computational
cost of  $\C{O}(r\, n\, d)$, and without changing the rank,
\begin{equation*}  
 \alpha\cdot \tnb{w} = \sum_{j=1}^{r} \alpha \bigotimes_{\nu=1}^d \vek{w}_{j}^{(\nu)}=
\sum_{j=1}^{r}  \bigotimes_{\nu=1}^d (\alpha_{\nu}\vek{w}_{j}^{(\nu)}) =
\sum_{j=1}^{r}  \bigotimes_{\nu=1}^d \vtil{w}_{j}^{(\nu)}, \quad
\text{with} \; \vtil{w}_{j}^{(\nu)}= {\alpha_{\nu}}\vek{w}_{j}^{(\nu)}.
\end{equation*}

\paragraph{The sum of two tensors} in the CP format 
$\tnb{w}=\tnb{u} + \tnb{v}$ can be written as follows
\begin{equation*}  
\tnb{w} = \tnb{u} + \tnb{v} =\left(\sum_{j=1}^{r_u} 
   \bigotimes_{\nu=1}^d \vek{u}_{j}^{(\nu)}\right) +
    \left(\sum_{k=1}^{r_v} \bigotimes_{\mu=1}^d \vek{v}_{k}^{(\mu)}\right)
  =
  \sum_{j=1}^{r_u+r_v} \bigotimes_{\nu=1}^d \vek{w}_{j}^{(\nu)} ,
\end{equation*}
where $\vek{w}_{j}^{(\nu)}:=\vek{u}_{j}^{(\nu)}$ for $j\leq r_u$ and 
$\vek{w}_{j}^{(\nu)}:=\vek{v}_{j}^{(\nu)}$ for $r_u< j\leq r_u+r_v$. 
The result generally has rank $r_u+r_v$, and as
this operation requires only concatenation of memory it  has
only a computing cost of $\C{O}(1)$.


\paragraph{The Hadamard product}
$\tnb{w}=\tnb{u}\odot \tnb{v}$ can be written as follows
\begin{equation*}
\tnb{w}=\tnb{u} \odot \tnb{v} =
  \left(\sum_{j=1}^{r_u} \bigotimes_{\nu=1}^d \vek{u}_{j}^{(\nu)}\right) 
\odot 
  \left(\sum_{k=1}^{r_v} \bigotimes_{\nu=1}^d \vek{v}_{k}^{(\nu)}\right)
=
\sum_{j=1}^{r_u}\sum_{k=1}^{r_v} \bigotimes_{\nu=1}^d \left(\vek{u}_{j}^{(\nu)} \odot
\vek{v}_{k}^{(\nu)}\right).
\end{equation*}
The new rank is generally $r_u\times r_v$, and the computational cost is 
$\C{O}(r_u \,r_v n\, d)$ arithmetic operations.

\paragraph{The Euclidean inner product}  is computed as follows:
\begin{equation*}
    \bkt{\tnb{u}}{\tnb{v}}_{\C{T}} =\bkt{\sum_{j=1}^{r_u} \bigotimes_{\nu=1}^d \vek{u}_{j}^{(\nu)}}%
{\sum_{k=1}^{r_v} \bigotimes_{\nu=1}^d \vek{v}_{k}^{(\nu)}}_{\C{T}}
=
\sum_{j=1}^{r_u}\sum_{k=1}^{r_v} \prod_{\nu=1}^d 
\bkt{\vek{u}_{j}^{(\nu)}}{\vek{v}_{k}^{(\nu)}}_{\C{P}_\nu}.
\end{equation*}
The computational cost of the inner product is $\C{O}(r_u\,r_v\,n\,d)$. 

We note that in all operations the numerical cost
grows only linearly with $d$, but the representation
rank of the resulting tensors may increase. 
Therefore, a rank truncation procedure is needed, which approximates a given 
tensor represented in the
canonical format with lower rank tensors up to a given accuracy.

\subsubsection{Rank truncation in the CP format}  \label{S:trunc}
Let $\tnb{w}$ be a tensor of rank $R$. 
Truncating $\tnb{w}$ to a new rank $r<R$ is a 
fundamental problem \citep{EspigDiss}.  This problem can be formulated as follows:
\begin{equation*}  
\text{find a}\; \tnb{w}^*\; \text{with rank}\;r \; \text{such that }
\forall \tnb{u} \;\;\text{with rank}\;r\; :\;\Vert 
\tnb{w} - \tnb{w}^*\Vert\leq \Vert \tnb{w} - \tnb{u} \Vert. 
\end{equation*}
Typical methods to solve this problem are the ALS-method and the Gauss-Newton-method.
The ALS method may show slow convergence, see \citep{EspigDiss, Espig2012VarCal} and references 
therein.  The Gauss-Newton method for $d\geq 3$ requires some additional 
assumptions, and also may not show any convergence at all \citep{EspigDiss}.

It is known, see for instance pp. 91--92 in \citep{khorBook18}, that the class of 
rank-$r$ CP tensors is a non-closed set in the corresponding tensor product space
for $d>2$.
Therefore there is no $\tnb{w}^*$ as above, and one may look for an
$\vepsilon$-solution, i.e.\ minimising within a deviation
of $\vepsilon$ from the infimum:
\begin{definition}[Approximation Problem]
\label{D:appTask2}
For a given tensor $\tnb{w}$ in CP format with rank $R$ and $\vepsilon >0$ we are looking for
minimal $r_\vepsilon < R$ and a tensor $\tnb{w}^*$ of rank $r_\vepsilon$
in CP format, such that:
\begin{equation*}
   \Vert \tnb{w} - \tnb{w}^* \Vert \leq \varepsilon \Vert \tnb{w} \Vert. 
\end{equation*}
\end{definition}
This problem is discussed in \citep{EspigDiss, ESHAROSCH09_1}. 

Properties of the TT format are described in Appendix~\ref{S:TT-rep}, and of the Tucker format in \cite{Espig2019}.

\section{Numerical examples} \label{S:num-exp}
%
  
The algorithms for post-processing high-dimensional data, or large volumes
of data are shown on a few illustrative examples.  Obviously, in any such
application the actual computing times and even the possibility of executing
such algorithms depend not only on that the huge amounts of data can be
compressed to reasonable size, but also on the possibility of executing
the algebraic operations with a reasonable speed on the compressed data.
Additionally, as the compression may deteriorate in the course of the computation,
one must assume that the intermediate results can be compressed to reasonable
size, and that this re-compression can be done efficiently.  All this depends
very much on the data, and currently there are only very few general results
to predict when this might or might not occur.  Some results in this direction
are already contained in \citep{ESHALIMA11_1, dolgov2014computation, DolgLitv15} 
for the low-rank TT-format which was also used
in our computations.

\subsection{Finding the value and location of the maximum element}  \label{SS:num-max}
Already in \reftab{tab:find_max} in \refS{intro} one could see an example
from \citep{EspigDiss}, illustrating the kinds of savings which are possible
when compared to the method of inspecting every element. 
We used low-rank TT-format storage, and the Algorithm~\ref{alg:power-RQ}.
In no case did it need more than 20 iterations.
The actual computing times on are shown in the fourth column of \reftab{tab:find_max},
and they behave like $\C{O}(n\,d)$.

\subsection{Finding level sets}  \label{SS:num-level}
We show two examples, the first is the same as in the previous \refSS{num-max}, wheres
the second one is a solution to a stochastic partial differential equation (SPDE).

\paragraph{Level set of the solution to a high-dimensional Poisson equation} is the
first example.  
\begin{table}[h!]
{\footnotesize
\caption{Computing times (4th column) on 3~GHz CPU to compute $\chi_S$.}
\label{tab:find_level}
\begin{center}
\begin{tabular}{r|l|l|r}
\hline
   & &                                       & actual 3~GHz CPU\\           
$d$ & \# elements.:  & $\approx$ years [a]   & iteration  time [s] for \\
    &   $N=n^d$   & inspect. $N$ & Newton-Schulz \\
\hline
 $25$  & $10^{50} $   & $10^{33} $  & $0.63$ \\
 $50$  & $10^{100}$   & $10^{83} $  & $1.30$ \\
$100$  & $10^{200}$   & $10^{133}$  & $2.66$ \\
$200$  & $10^{400}$   & $10^{233}$  & $6.03$ \\
$400$  & $10^{800}$   & $10^{433}$  & $12.97$ \\
$800$  & $10^{1600}$  & $10^{1233}$ & $25.32$ \\
\hline
\end{tabular}
\end{center}
}
\end{table}%
The data vector is the same one as used in \refS{intro}.  
But now we want to find the indices of the elements which are in the set 
$S=\,]\!-\infty,0.25[$, and according to \refSS{post-proc-tasks}
one has to compute the characteristic function $\chi_S$ as given in \refeq{eq:char-def}.

For this one needs the sign function given in \refeq{eq:sign-def}, computed with
the Newton-Schulz iteration function in \refeq{eq:iter-sign-N-S}.  The results are
shown in \reftab{tab:find_level}.  The first three columns are similar to those
in \reftab{tab:find_max}, and the fourth column are the actual computing times in
seconds on a 3~GHz CPU, which behave like $\C{O}(n\,d)$.

\paragraph{Level sets of the solution to a stochastic partial differential equation} (SPDE)
is the next example.  This example is described in a bit more detail, see also
\citep{ESHALIMA11_1}.

Two scientific software libraries \textit{sglib} and \textit{TensorCalculus} 
\citep{sglib, TensorCalEspig} were used.  With the \textit{sglib} procedures we discretise 
the SPDE \refeq{eq:model-problem}, and with the \textit{TensorCalculus}
we solve the obtained tensor equation and compute the quantities of interests.

Consider the following diffusion equation with uncertain diffusion coefficient $\kappa(\omega,x)$:
\begin{equation}
\left.\begin{alignedat}{2}-\nabla(\kappa(\omega,x)\nabla u(\omega,x)) & = f(\omega,x) & 
\quad & \text{ in } \C{G},\\
u(\omega,x) & =0 &  & \text{ on } \partial\C{G},
\end{alignedat}
\ \right\} \ \text{a.s. in \ensuremath{\omega\in\Omega}},\label{eq:model-problem}
\end{equation}
Here $\C{G}$ is the two-dimensional L-shaped domain $[-1,1]^2 \setminus [0,1]^2$,
and a triangular mesh with $557$ mesh points was used for the spatial piece-wise
linear finite element discretisation, and including the boundary values,
the discrete spatial solution needs $M_d = 557$ storage locations.

After a stochastic discretisation with the \KL{} and Polynomial Chaos Expansion,
and applying the stochastic Galerkin method as in \citep{Matthies_Keese_2005CMA_SFEM,
Matthies-2008-ZAMM, dolgov2014computation, DolgLitv15, ESHALIMA11_1, MeWhAlHgmPw12-p}, 
we can compute the solution $u(\omega,x)$ in a chosen tensor representation.

The random  field $\kappa(\omega,x)$ was taken to have a
shifted log-normal distribution for, i.e.\ 
$\log(\kappa(\omega,x)-1.1)$ has a normal distribution with parameters 
$\{\mu=0.5,\, \sigma^2=1.0\}$.  The isotropic and homogeneous covariance function is of 
Gaussian type with covariance lengths $\ell_x=\ell_y=0.3$. 

Ten \KL{} terms were taken to represent the mean zero random part of 
the field $\kappa(\omega,x)$, 
and for the polynomial chaos expansion (PCE) multivariate 
Hermite polynomials of maximum second degree were taken.

For the right-hand side $f(\omega,x)$ a $\beta$-distribution $\{4,2\}$  was chosen
for the random field.  The covariance function is also of 
Gaussian type with covariance lengths $\ell_x=\ell_y=0.6$.
Again ten terms in the \KL{} expansion for the zero mean random part of
$f(\omega,x)$ with a maximum
second degree polynomial chaos were taken.

The total stochastic dimension of the
solution $u(\omega,x)$ is $20$ --- ten random variables from the ten terms of the
\KL{} expansion of the coefficient $\kappa(\omega,x)$, plus ten terms of the
\KL{} expansion of the right hand side $f(\omega,x)$ ---
i.e.\ the multi-index $\alpha$
will consist of $20$ indices $(\alpha=(\alpha_1,...,\alpha_{20}))$, and the
solution will be represented as a function of $20$ independent random variables.
Together with the extra dimension from the spatial discretisation --- which is lumped
into one --- one has a discrete version of a function $u(\omega,x)$
 on $d=21$ dimensions.  With second degree polynomials, there are $M_s =3$ polynomials
 for each random variable, thus the stochastic number degrees of freedom
 is  $N_s = M_s^{20} =  3^{20} = 3,486,784,401 \approx 3.5\times 10^9$.
Thus the total number of entries for the full tensor is $N_d \times N_s = N_d \times M_s^{20}
= 557 \times 3,486,784,401 = 1,942,138,911,357 \approx 2\times 10^{12}$,
and the discrete solution  $\tnb{u} \in \R^{557} \ten \Ten_{\mu=1}^{20} \R^{3}
\cong \D{R}^{1,942,138,911,357}$.  Thus, even for this relatively coarse
discretisation, full storage of just the solution vector
would require approximately 16~TB of memory.

For the numerical solution, a low-rank tensor solver was used \citep{Zander10, 
dolgov2014computation, DolgLitv15, ESHALIMA11_1, MeWhAlHgmPw12-p}, and the solution was
represented in the CP-format \refeq{eq:CP-format} with a rank of $r=231$.
Thus
\begin{equation*}
\tnb{u}=\sum_{j=1}^{231} \vek{u}_{j0} \otimes \Ten_{\mu=1}^{20} \vek{u}_{j \mu} \in \R^{557}
\ten \Ten_{\mu=1}^{20} \R^{3},
\end{equation*}
and the storage of the solution vector requires only $231\times(557 + 20\times 3)=142,527$ 
storage locations, i.e.\ less than 1.15~MB.
\begin{table}[h!]
{\footnotesize
\caption{Computing $\chi_S(\tnb{u})$.}
\label{tab:find_level-SPDE}
\begin{center}
\begin{tabular}{c|c|c|c|l}
  $b/\|\mat{u}\|_\infty$   &   rank $\chi_S$ & max it rank  & $k_{\max}$ its & error\\
\hline
 $0.2$  & $ 12 $  & $ 24 $ & $ 12 $ & $ 2.9\times 10^{-8}$\\
 $0.4$  & $ 12 $  & $ 20 $ & $ 20 $ & $ 1.9\times 10^{-7}$\\
 $0.6$  & $  \phantom{0}8 $  & $ 16 $ & $ 12 $ & $ 1.6\times 10^{-7}$\\
 $0.8$  & $  \phantom{0}8 $  & $ 15 $ & $ \phantom{0}8 $ & $ 1.2\times 10^{-7}$\\
\hline
\end{tabular}
\end{center}
}
\end{table}%

Finally, we computed $\chi_S(\tnb{u})$ for $S=]-\infty,b[$ by computing
 $\sign(b\|\tnb{u}\|_\infty \tnb{1} -\tnb{u})$ for $b \in \{0.2,\, 0.4,\, 0.6,\,0.8\}$,
which is shown in the first column of Table~\ref{tab:find_level-SPDE}.
The final representation ranks of $\chi_S(\tnb{u})$ resp.\ 
$\sign(b \nd{\tnb{u}}_\infty \tnb{1} -\tnb{u})$ are
given in the second column.  In this numerical example, the ranks are in all cases
smaller then $13$.  The $\sign$-function was computed with the iteration from 
Algorithm~\ref{alg:basic}.  The iteration function is given in \refeq{eq:iter-sign-N-S},
and the iteration terminated after $k_{\max}$ steps, shown in column four.
The maximal representation rank of the iterates is documented in the third column.  
The error $\nd{\tnb{1} - \tnb{u}_{k_{\max}}\odot 
\tnb{u}_{k_{\max}}}/\nd{(b \nd{\tnb{u}}_\infty \tnb{1} - \tnb{u})}$ is given in the
last column.  Each computation --- each row of the table --- took less than $10$ min
on a 3~GHz CPU.

%
%
%
%
%
%
%
%
%
%
%
%




\section{Conclusion} \label{S:concl}
Large volumes of data are used resp.\ generated in an increasing number of
instances in science, technology, and business; some of which were
mentioned in \refSS{xmpl-motiv}.  Often this data is in the form of
a sample of a high dimensional function, or it can be reshaped in this way.
It may well happen that it is not possible to store such data in its entirety,
so that it has to be compressed in some way.  And it is usually not enough to 
compress these large data sets, efficient numerical 
algorithms are required to post-process such data in a compressed data format.
All these post-processing tasks may be trivial in low dimensions and for small
data sets, but they become non-trivial for high-dimensional data, 
say $d>5$, which is in some compressed format. 
We have introduced some numerical methods which may partially close this gap.

In this work we assumed that 1) the data are given or can be approximated in a compressed data format (e.g., in a low-rank tensor format); 2) a recompression (tensor rank truncation in our case) procedure is available and the ranks after linear algebra operations are not increases too much.

Here we have formulated algorithms for such post-processing, which only use
the structure of an abstract commutative algebra.  The Hadamard algebra
is an example of such an algebra, and it can be defined for all real valued
data sets.  The Hadamard product of two tensors turned out to be a key tool in 
approximating the variance, point-wise inverse, level sets, frequency, and the 
maximal element of huge multi-dimensional data sets.

In our example application, the compression is low-rank tensor
approximation, and we showed how to perform the algebra operations in
some example low-rank tensor formats.
We assumed that the initial data set has low-rank 
tensor representation, and showed that this low-rank property
can be preserved during the whole computing process. 

The choice of the tensor format is not crucial; it could be, for example, the
canonical polyadic format, the Tucker format, the tensor train format, 
and many others.  The only requirement is that it should be possible to 
compute the operations of the Hadamard algebra in a reasonable (linear in $n$) time,
and the existence of stable rank truncation algorithms to truncate 
possibly large intermediate ranks.
The algorithms needed were among others the multiplicative inverse and the sign
function, and this could be fashioned after the algorithms for matrix algebras.
Finally, the algorithms are demonstrated on some high dimensional data which comes
from the solution of high dimensional or parametric resp.\ stochastic elliptic
partial differential equations.

\appendix

\section{The Tensor Train format} \label{S:TT-rep}
The tensor train (TT) format is described in \citep{OSTY09, OS11, HA12,  khorBook18}.
As already noted, it was originally developed in quantum chemistry as ``matrix product
states'' (MPS), see  \citep{VI03} and references therein, and rediscovered later
\cite{oseledetsTyrt2010, oseledets2011}.
For a motivation, let us start with a well-known example \citep{Kazeev_Laplace_TT}.
\begin{example}
Consider the $d$-dimensional Laplacian operator discretised 
with standard finite differences over a uniform 
tensor grid with $n$ degrees of freedom in each direction.
It has the Kronecker (canonical) rank-$d$ representation:
\begin{equation*}
\tnb{A} = \vek{A}\otimes 
\overbrace{\vek{I} \otimes \cdots \otimes \vek{I}}^{d-1 \text{\rm{} \; times}} + 
\vek{I} \otimes \vek{A} \otimes \vek{I} \cdots \otimes \vek{I}  
+ \cdots + \vek{I}\otimes\cdots  \otimes \vek{I}\otimes \vek{A} \in \D{R}^{n^d \times n^d}
\end{equation*}
with $\vek{A}= \text{\rm tridiag}\{-1,2,-1\}\in \D{R}^{n \times n}$, and $\vek{I}$ being the 
$n \times n$ identity.
However, the same operator in the TT format is explicitly representable with 
all TT ranks equal to 2 for any dimension \citep{Kazeev_Laplace_TT},
\begin{equation*}
\tnb{A}=  
 \begin{pmatrix}
  \vek{A}, & \vek{I}
 \end{pmatrix}
 \Join
 \overbrace{
 \begin{pmatrix}
  \vek{I} & \vek{0} \\
  \vek{A} & \vek{I}
 \end{pmatrix}\Join
\dots
 \Join
 \begin{pmatrix}
  \vek{I} & \vek{0} \\
  \vek{A} & \vek{I}
 \end{pmatrix}
 }^{d-2 \text{\rm{} \; times}}
 \Join
 \begin{pmatrix}
  \vek{I}  \\
  \vek{A}
 \end{pmatrix},
 \end{equation*}
where the ``strong Kronecker'' product operation `` $\Join $ '' is defined as a regular matrix 
product for the first level of TT cores, i.e.\ the block matrices, 
and the inner blocks  are multiplied by means of the Kronecker or tensor product;
for example 
\[
 \begin{pmatrix}  \vek{A}, & \vek{I} \end{pmatrix}  \Join  \begin{pmatrix}
  \vek{I} & \vek{0} \\   \vek{A} & \vek{I} \end{pmatrix} = \begin{pmatrix}  
  \vek{A} \otimes \vek{I}+\vek{I} \otimes \vek{A}, & \vek{I}\otimes \vek{I} \end{pmatrix} .
\]
In the element-wise matrix product notation, the Laplace operator reads
\[
\tnb{A}[\vek{i},\vek{j}] =
\begin{pmatrix}\vek{A}[i_1,j_1] & \vek{I}[i_1,j_1]\end{pmatrix}
\begin{pmatrix}
 \vek{I}[i_2,j_2] & 0 \\
 \vek{A}[i_2,j_2] & \vek{I}[i_2,j_2]
 \end{pmatrix} \cdots
 \begin{pmatrix}
  \vek{I}[i_d,j_d]  \\
  \vek{A}[i_d,j_d]
 \end{pmatrix}.
\]
\end{example}
%
Thus, the $d$-dimensional Laplacian operator is separable with TT ranks equal to 2,
and each TT component (core) is defined by a one-dimensional Laplacian.

\begin{definition}[TT-Format, TT-Representation, TT-Ranks]\label{D:TTFormat}
The \emph{TT-tensor format} is for variable \emph{TT-representation
ranks} $\vek{r}=(r_0, \dots, r_{d})\in \N^{d+1}$ --- with $r_0 = r_d = 1$ and
under the assumption that $d>2$ --- defined by the
following multilinear mapping
\begin{align}   \label{eq:TTRepresentation}
  \tns{U}_{\TT} : \C{P}_{\TT,\vek{r}} &:=  \bigtimes_{\nu=1}^{d} 
  \C{P}_\nu^{r_{\nu-1} \times r_{\nu}}
  \rightarrow   \C{T},\quad \C{P}_\nu = \D{R}^{M_\nu}\; (\nu=1,\dots,d),   \\ \nonumber
   \C{P}_{\TT,\vek{r}}\ni\tnb{P} &=(\tnb{W}^{(\nu)} = (\vek{w}^{(\nu)}_{j_{\nu-1} j_\nu})
   \in\C{P}_\nu^{r_{\nu-1} \times r_{\nu}} :1\le j_\nu\le r_\nu, 1\le\nu\le d) \\ \nonumber
   &\qquad\mapsto 
   \tns{U}_{\TT}(\tnb{P}):=\tnb{w}=\sum_{j_0=1}^{r_0}\dots\sum_{j_d=1}^{r_d} 
     \Ten_{\nu=1}^{d} \vek{w}_{j_{\nu-1} j_\nu}^{(\nu)} \; \in \C{T}.
\end{align}
We call  $\tnb{w}:=(\tns{w}_{i_1\dots i_d})=U_{\TT,\vek{r}}(\tnb{P})$ a tensor 
represented in the train tensor format.
Note that the TT-cores $\tnb{W}^{(\nu)}$ may be viewed as a vector 
valued $r_{\nu-1}\times r_\nu$ matrix
with the vector $\vek{w}^{(\nu)}_{j_{\nu-1} j_\nu}\in\D{R}^{M_\nu}$ with the components
$w^{(\nu)}_{j_{\nu-1} j_\nu}[i_\nu]: 1\le i_\nu \le M_\nu$
at index position $(j_{\nu-1}, j_\nu)$.    The representation in components is then
\begin{equation}   \label{eq:TTRep-comp-1}
  (\tns{w}_{i_1\dots i_d})=\sum_{j_0=1}^{r_0} \dots \sum_{j_{d}=1}^{r_{d}} 
  w_{j_0 j_1}^{(1)}[i_1] \cdots 
  w_{j_{\nu-1} j_\nu}^{(\nu)}[i_{\nu}] \cdots w_{j_{d-1} j_d}^{(d)}[i_d]
\end{equation}
Alternatively, each TT-core $\tnb{W}^{(\nu)}$ may be seen as a vector
of $r_{\nu-1}\times r_\nu$ matrices $\vek{W}^{(\nu)}_{i_\nu}$ of length $M_\nu$, i.e.\
$\tnb{W}^{(\nu)} = (\vek{W}^{(\nu)}_{i_\nu}):1\le i_\nu \le M_\nu$.  Then
the representation \refeq{eq:TTRep-comp-1} reads
\begin{equation}   \label{eq:TTRep-comp-2}
  (\tns{w}_{i_1\dots i_d})= \prod_{\nu=1}^d \vek{W}^{(\nu)}_{i_\nu},
\end{equation}
which explains the name \emph{matrix product state}.  Observe that the first matrix
is always a row vector as $r_0=1$, and the last matrix is always a column vector as $r_d=1$.
The matrix components $\vek{W}^{(\nu)}_{i_\nu}$ of the TT-cores are also called ``carriages''
or ``waggons'' with ``wheels'' $i_\nu$ at the bottom, coupled to the next ``carriage''
or ``waggon'' via the matrix product.  This explains the \emph{tensor train} name.
If one notes more carefully $\tnb{W}^{(\nu)} \in \C{P}_\nu^{r_{\nu-1} \times r_{\nu}} = 
\D{R}^{r_{\nu-1}} \otimes \D{R}^{M_\nu} \otimes \D{R}^{r_\nu}$, then \refeq{eq:TTRep-comp-2}
can be written more concisely as
\begin{equation}   \label{eq:TTRep-comp-3}
  \tnb{w} = \tns{U}_{\TT}(\tnb{P}) =  \tnb{W}^{(1)} \times_3^1 \tnb{W}^{(2)} \times_3^1 \cdots
  \times_3^1 \tnb{W}^{(d)} ,
\end{equation}
where $\tnb{U} \times_k^\ell \tnb{V}$ is a contraction of the $k$-th index of $\tnb{U}$
with the $\ell$-th index of $\tnb{V}$, where one often writes just $\times_k$ for
$\times^1_k$.  Thus in \refeq{eq:TTRep-comp-3} the contractions leave the indices from
the $\D{R}^{M_\nu}$ untouched, so that the tensor $\tnb{w}$ is formed.
\end{definition}

Each TT-core (or \emph{block}) $\tnb{W}^{(\nu)}$ is defined by 
$r_{\nu-1}\times r_{\nu} \times M_\nu$ numbers. Assuming $n=M_\nu$ for all $\nu=1,\ldots,d$,
the total number of entries scales as $\mathcal{O}(d\,n\,r^2)$, which is tractable 
as long as $r=\max\{r_k\}$ is moderate.

\begin{figure}[h]
\centering
\includegraphics[width=0.7\textwidth]{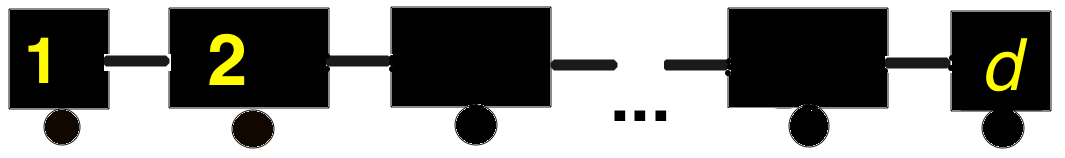}
\caption{Schema of the TT tensor decomposition.
The waggons denote the TT cores and each wheel denotes
the index $i_\nu$. Each waggon is connected with neighbours by indices $j_{\nu-1}$ and $j_{\nu}$.}
\label{fig:TT}
\end{figure}

A pictorial representation of the schema for the TT tensor format is shown in Figure~\ref{fig:TT}.
It shows $d$ connected waggons with one wheel.  The waggons denote the TT-cores,
and each wheel denotes the index $i_\nu$.  The waggons for $\nu = 2,\dots,(d-1)$ are connected
with their neighbours by two indices $j_{\nu-1}$ and $j_{\nu}$.  The first and the last
waggons are connected by only one index, namely $j_1$ and $j_{d-1}$ respectively.  Since by 
the convention in \refD{TTFormat} above, $r_0=r_d=1$, the indices $j_0$ and $j_d$ run from 1 to 1, i.e.\ are purely formal.

\ignore{
All the TT-cores $\tnb{W}^{(\nu)}$ are tensors of degree 3 --- i.e.\
$\vek{W}^{(\nu)}_{i_\nu}[j_{\nu-1}, j_\nu]$ --- with indices
$j_{\nu-1}$, $i_\nu$, and $j_\nu$, representing the carriages or waggons.
The $i_\nu$ index is the wheel, and the $j_{\nu-1}$ and $j_\nu$ indices are the
``hitch'' to the preceding and the following waggon or carriage.

The first waggon in Figure~\ref{fig:TT} symbolises the TT-core $\vek{W}^{(1)}_{i_1}$,
and as $1 \le j_0 \le r_0=1$, it is only formally a tensor of degree 3, but rather
one of degree 2, i.e.\ a matrix, as the $j_0$ index is constant and hence does not really count.
The second waggon is the TT-core $\vek{W}^{(2)}_{i_2}$, and the last waggon is $\vek{W}^{(d)}_{i_d}$. 
The first waggon is connected with the second by the ``hitch'', the index $j_1$, the second with the 
third one by the index $j_2$, and the penultimate one to the last by the index $j_{d-1}$, 
see \refeq{eq:TTRep-comp-1} and \refeq{eq:TTRep-comp-3}.  Again, as $1 \le j_d \le r_d=1$,
the last waggon is really a matrix.  
}
The waggon or carriage 2 --- $\tnb{W}^{(2)}$ ---
is a tensor of degree 3, described by three 
indices $j_1$ (the left hitch), $i_2$ (the wheel), and $j_2$ (the right hitch).  
Multiplication of the third TT-core with the second 
and forth cores means the tensor contraction by the indices $j_1$ and $j_2$.  If we perform 
tensor contraction of all TT-cores over the indices $j_1,\ldots,j_{d-1}$,
and disregard the purely formal constant indices
$j_0$ and $j_d$, then the indices $j_0,\ldots,j_{d}$ --- the hitches --- will disappear, and 
only the indices $i_1,\ldots,i_d$ --- the wheels --- will be left.



\subsection{Basic operations with the TT format}
We follow to the work of Oseledets \citep{OS11} and list the major properties of the 
TT-tensor format.

\paragraph{The multiplication with scalar} $\alpha$ could be simply
done by multiplying one of the TT-cores $\tnb{W}^{(\nu)}$ in the
representation \refeq{eq:TTRep-comp-3} for any $\nu$ in
$\tnb{w} =  \tnb{W}^{(1)} \times_3^1 \tnb{W}^{(2)} \times_3^1 \cdots
  \times_3^1 \tnb{W}^{(d)}$.
But to balance the effect better, define $\alpha_{\nu} := \sqrt[d]{|\alpha|}$ for all $\nu > 1$,
and $\alpha_{1} := \sign(\alpha)\sqrt[d]{|\alpha|}$.  Then $\Ttil{w} = \alpha\cdot\tnb{w}$
is given by
\[  \Ttil{w} =  (\alpha_1\cdot\tnb{W}^{(1)}) \times_3^1 
  (\alpha_2\cdot\tnb{W}^{(2)}) \times_3^1 \cdots  \times_3^1 (\alpha_d\cdot\tnb{W}^{(d)})
  = \Ttil{W}^{(1)}\times_3^1 \cdots \times_3^1 \Ttil{W}^{(d)}.
\]
The new cores are given by $\Ttil{W}^{(\nu)} = (\vtil{W}^{(\nu)}_{i_\nu})
= (\alpha_\nu \vek{W}^{(\nu)}_{i_\nu})$, a sequence of new ``carriage'' matrices. 
The computational complexity is $\mathcal{O}(d\,n\,r^2)$.

\paragraph{Addition of two TT-tensors}
Assume two tensors $\tnb{u}$ and $\tnb{v}$ are given in the TT-tensor format
as in \refeq{eq:TTRep-comp-2}, i.e.\
$(\tns{u}_{i_1\dots i_d})= \prod_{\nu=1}^d \vek{U}^{(\nu)}_{i_\nu}$ and 
$(\tns{v}_{i_1\dots i_d})= \prod_{\nu=1}^d \vek{V}^{(\nu)}_{i_\nu}$.  The sum
$\tnb{w}=\tnb{u}+\tnb{v}$ is given by the new cores $\vek{W}^{(\nu)}_{i_\nu}$ such that
$(\tns{w}_{i_1\dots i_d})= \prod_{\nu=1}^d \vek{W}^{(\nu)}_{i_\nu}$, where
\begin{equation*}
    \vek{W}^{(\nu)}_{i_\nu}=
\begin{pmatrix}
\vek{U}^{(\nu)}_{i_\nu} & \vek{0} \\
 \vek{0} & \vek{V}^{(\nu)}_{i_\nu} 
\end{pmatrix}, \quad 1\le i_\nu \le r_\nu, 2\le\nu\le d-1.
\end{equation*}
and the first and the last cores will be 
\[
\vek{W}^{(1)}_{i_1}=\left(\vek{U}^{(1)}_{i_1}\; \vek{V}^{(1)}_{i_1} \right)  \quad \text{ and } 
\quad
\vek{W}^{(d)}_{i_d}=\left( \begin{array}{c} \vek{U}^{(d)}_{i_d} \\ \vek{V}^{(d)}_{i_d}
\end{array} \right).
\]
As only storage may have to be concatenated, the computational cost is $\C{O}(1)$,
but as  the carriages resp.\ TT-cores grow, the final rank will generally be the
sum of the ranks.

\paragraph{The Hadamard product} $\tnb{w}=\tnb{u} \odot \tnb{v}$ in the TT format
 is computed as follows.  Assume two tensors $\tnb{u}$ and $\tnb{v}$ are given in 
 the TT tensor format as in \refeq{eq:TTRep-comp-2}, i.e.\
$(\tns{u}_{i_1\dots i_d})= \prod_{\nu=1}^d \vek{U}^{(\nu)}_{i_\nu}$ and 
$(\tns{v}_{i_1\dots i_d})= \prod_{\nu=1}^d \vek{V}^{(\nu)}_{i_\nu}$.  The Hadamard product is
\[
    (\tns{w}_{i_1\dots i_d}) = (\tns{u}_{i_1\dots i_d} \cdot \tns{v}_{i_1\dots i_d}) .
\]
The tensor $\tnb{w}$ has also the TT-tensor format, namely with the new cores
\[
   \vek{W}^{(\nu)}_{i_\nu} = \vek{U}^{(\nu)}_{i_\nu}\otimes_{\mrm{K}} \vek{V}^{(\nu)}_{i_\nu}, 
   \quad 1\le i_\nu \le r_\nu, 1\le\nu\le d ,
\]
where $\otimes_{\mrm{K}}$ is the Kronecker product of two matrices \citep{HA12}.
The rank of $\tnb{W}^{(\nu)} = (\vek{W}^{(\nu)}_{i_\nu})$ is the product
of the ranks of the TT-cores  $\tnb{U}^{(\nu)}$ and $ \tnb{V}^{(\nu)}$.

\paragraph{The Euclidean inner product} of two tensors in the TT-format
as in \refeq{eq:TTRepresentation} 
\begin{equation*}
\tnb{u} = \sum_{j_0=1}^{r^u_0}\dots\sum_{j_d=1}^{r^u_d} 
     \Ten_{\nu=1}^{d} \vek{u}_{j_{\nu-1} j_\nu}^{(\nu)} , \qquad
\tnb{v} = \sum_{j_0=1}^{r^v_0}\dots\sum_{j_d=1}^{r^v_d} 
     \Ten_{\nu=1}^{d} \vek{v}_{j_{\nu-1} j_\nu}^{(\nu)} ,
\end{equation*}
with ranks $\vek{r}^u$ and $\vek{r}^v$ can be computed as follows:
\[
\bkt{\tnb{u}}{\tnb{v}}_{\C{T}}= \sum_{j_0=1}^{r^u_0}\dots\sum_{j_d=1}^{r^u_d}
  \sum_{i_0=1}^{r^v_0}\dots\sum_{i_d=1}^{r^v_d}  \prod_{\nu=1}^d \;
       \bkt{\vek{u}_{j_{\nu-1} j_\nu}^{(\nu)}}{\vek{v}_{i_{\nu-1} i_\nu}^{(\nu)}}_{\C{P}_\nu}.
\]
The computational complexity is $\C{O}(d\,n\,r^4)$, and can be reduced further \citep{OselTT}. 
%
%
%

\subsection{Rank truncation in the TT format}
The rank truncation opearation is based on the SVD algorithm and requires 
$\C{O}(d\,n\,r^3)$ operations \citep{GrasHack11}.  The TT-rounding algorithm (p. 2305 in \citep{oseledets2011}) is based on QR decomposition and costs $\C{O}(d\,n\,r^3)$.

Corollary 2.4 in \citep{oseledets2011} states that for a given tensor $\tnb{w}$ and
rank bounds $r_k$, the best approximation to $\tnb{w}$ in the Frobenius norm with TT-ranks 
bounded by $r_k$ always exist (denote it by $\tnb{w}^*$), and the TT-approximation 
$\tnb{u}$ computed by the TT-SVD algorithm (p. 2301 in \citep{oseledets2011}) is quasi-optimal:
\begin{equation}
    \Vert \tnb{w}-\tnb{u} \Vert_F\leq \sqrt{d-1} \Vert \tnb{w} - \tnb{w}^{*} \Vert_F.
\end{equation}
In \cite{Kressner17_Trunc} the authors suggested a new re-compression 
randomised algorithm for Tucker and TT tensor formats.

\cleardoublepage 

\begin{footnotesize}
\providecommand{\bysame}{\leavevmode\hbox to3em{\hrulefill}\thinspace}
\providecommand{\MR}{\relax\ifhmode\unskip\space\fi MR }
\providecommand{\MRhref}[2]{%
  \href{http://www.ams.org/mathscinet-getitem?mr=#1}{#2}
}
\providecommand{\href}[2]{#2}

\end{footnotesize}





\end{document}